\documentclass[12pt]{amsart}

\usepackage{amsmath} 
\usepackage{amssymb}

\makeatletter
\@namedef{subjclassname@2020}{%
  \textup{2020} Mathematics Subject Classification}

\newtheorem{theorem}{Theorem}[section] 
\newtheorem{claim}{Claim}[theorem]
\newtheorem{lemma}[theorem]{Lemma} 
\newtheorem{proposition}[theorem]{Proposition} 
\newtheorem{corollary}[theorem]{Corollary} 
\newtheorem{observation}[theorem]{Observation} 

\theoremstyle{definition}
\newtheorem{problem}[theorem]{Problem} 
 
\newtheorem{definition}[theorem]{Definition}
\newtheorem{hypothesis}[theorem]{Assumptions}
\newtheorem{example}[theorem]{Example}

\theoremstyle{remark}

\newcount\skewfactor
\def\mathunderaccent#1#2 {\let\theaccent#1\skewfactor#2
\mathpalette\putaccentunder}
\def\putaccentunder#1#2{\oalign{$#1#2$\crcr\hidewidth
\vbox to.2ex{\hbox{$#1\skew\skewfactor\theaccent{}$}\vss}\hidewidth}}
\def\name{\mathunderaccent\tilde-3 }

\newcommand{\forces}{\Vdash} 
\newcommand{\con}{{\mathfrak c}} 
\newcommand{\can}{{}^{\omega}2} 
 
\newcommand{\rest}{{\restriction}}
 
\newcommand{\rng}{{\rm rng}}

\newcommand{\vtl}{\vartriangleleft} 
\newcommand{\vare}{\varepsilon} 
\newcommand{\fstwo}{{}^{\omega>}2}
\newcommand{\per}{{\rm per}}
\newcommand{\cf}{{\rm cf}}
\newcommand{\rk}{{\rm rk}}
\newcommand{\stnd}{{\rm stnd}}
\newcommand{\ndrk}{{\rm ndrk}}
\newcommand{\NDRK}{{\rm NDRK}}
\newcommand{\Mtk}{{{\mathbf M}_{\bar{T},\bar{\cO}}}}
\newcommand{\fMtk}{{{\mathbf M}^n_{\bar{t},\bar{\cO}^6}}}

\newcommand{\cA}{{\mathcal A}}

\newcommand{\cB}{{\mathcal B}}

\newcommand{\bbC}{{\mathbb C}}
\newcommand{\cC}{{\mathcal C}}

\newcommand{\gd}{{\mathfrak d}}

\newcommand{\cL}{{\mathcal L}}

\newcommand{\cM}{{\mathcal M}}
\newcommand{\bM}{{\mathbf M}}
\newcommand{\bbM}{{\mathbb M}}
\newcommand{\bm}{{\mathbf m}}
\newcommand{\bn}{{\mathbf n}}
\newcommand{\cO}{{\mathcal O}}

\newcommand{\bbP}{{\mathbb P}}
\newcommand{\bbQ}{{\mathbb Q}}

\newcommand{\bV}{{\mathbf V}}

\newcommand{\bbZ}{{\mathbb Z}}

\begin{document}

\title[Overlapping translations IV]{Borel sets without perfectly many
  overlapping translations IV} 

\author{Andrzej Ros{\l}anowski}
\address{Department of Mathematics\\
University of Nebraska at Omaha\\
Omaha, NE 68182-0243, USA}
\email{aroslanowski@unomaha.edu}

\author{Saharon Shelah}
\address{Institute of Mathematics\\
 The Hebrew University of Jerusalem\\
 91904 Jerusalem, Israel\\
 and  Department of Mathematics\\
 Rutgers University\\
 New Brunswick, NJ 08854, USA}
\email{shelah@math.huji.ac.il}
\urladdr{http://shelah.logic.at}
\thanks{Publication 1240 of the second author. Research partially supported
  by the Israel Science Foundation (ISF) grant no: 1838/19\\ 
  Both authors are grateful to an individual who prefers to remain anonymous
  for providing typing  services that were used during the work on the paper. }

 
\begin{abstract}
We show that, consistently, there exists a  Borel set $B\subseteq\can$
admitting a sequence $\langle\eta_\alpha:\alpha<\lambda\rangle$ of distinct
elements of $\can$ such that $(\eta_\alpha+B)\cap (\eta_\beta+B)$ is
uncountable for all $\alpha,\beta<\lambda$ but with no perfect set $P$ such
that $|(\eta+B)\cap (\nu+B)|\geq 6$ for any distinct $\eta,\nu\in P$.
This answers two questions from our previous works, \cite[Problem
5.1]{RoSh:1138}, \cite[Problem 7.6]{RoSh:1187}.
\end{abstract}



\maketitle

\section{Introduction}
In the series of articles \cite{RoSh:1138, RoSh:1170, RoSh:1187} we
investigated the existence of Borel sets with many, but not too many pairwise
non-disjoint translations. For instance, in \cite{RoSh:1170}, for a
countable ordinal $\vare<\omega_1$ and an integer $2\leq \iota<\omega$ we
constructed a $\Sigma^0_2$ set $B\subseteq \can$ with the following
property. 
\begin{quotation}
{\em In some ccc forcing notion there is a sequence
$\langle\rho_\alpha:\alpha<  \aleph_\vare\rangle$ of distinct elements of
$\can$ such that 
\[\big|(\rho_\alpha+B)\cap (\rho_\beta+B)\big|\geq 2\iota\mbox{ for all  
    }\alpha,\beta<\lambda\]
but in no extension there is a perfect set of such $\rho$'s. }
\end{quotation}
Similar resuts for the general case of perfect Abelian Polish groups were
presented in \cite{RoSh:1187}. However, in all those cases when discussing
nonempty intersections we considered finite intersections only. It seemed
that our arguments really needed a finite enumeration of ``witnesses for 
nondisjointness''. So in   \cite[Problem 5.1]{RoSh:1138} and  \cite[Problem
7.6]{RoSh:1187} we asked if there is a ccc forcing notion $\bbP$ adding a
$\Sigma^0_2$ subset $B$ of the Cantor space $\can$ such that     
\begin{quotation}
{\em for some $H\subseteq \can$ of size $\lambda$, the intersections 
  $(B+h)\cap (B+h')$ are infinite (uncountable, respectively) for all
  $h,h'\in H$, but for every perfect set $P\subseteq \can$ there are
  $x,x'\in P$ with the intersection $(B+x)\cap (B+x')$ finite (countable,
  respectively).} 
\end{quotation}
In the present paper we answer the above two questions positively. Our
forcing construction slightly generalizes and simplifies that of
\cite{RoSh:1138, RoSh:1170}. This allows us to show a stronger result:  
\begin{quotation}
{\em If $\lambda<\lambda_{\omega_1}$ then some ccc forcing notion adds a
  $\Sigma^0_2$ set $B$ which has $\lambda$ translations with pairwise
  uncountable intersections, while for every perfect set $P\subseteq \can$
  there are $x,x'\in P$ with  $|(B+x)\cap (B+x')|<6$.} 
\end{quotation}
The article is organized as follows. First, in Section 2, we recall the
splitting rank from Shelah \cite{Sh:522}. This rank was fundamental for the
question of no perfect squares and it is fundamental for problems of 
nondisjoint translations as well. Then, in the third section we introduce
nice indexed bases $\bar{\cO}$ and we define when translations of a
$\Sigma^0_2$ set have $\bar{\cO}$--large intersection. This allows us to put
in the same framework sets with finite, infinite and uncountable
intersections. We also analyze when a $\Sigma^0_2$ set may have a perfect
set of translations with $\bar{\cO}$--large intersections and we introduce a
non-disjointness rank on finite approximations. Our main consistency theorem
is presented in the fourth section. In the final part of the paper we
summarize our results and pose a few relevant problems. 
\medskip

\noindent{\bf Notation}: Our notation is standard and
compatible with that of classical textbooks (like Jech \cite{J} or
Bartoszy\'nski and Judah \cite{BaJu95}). However, in forcing we keep the
older convention that {\em a stronger condition is the larger one}.

\begin{enumerate}
\item For a set $u$ we let $u^{\langle 2\rangle}=\{(x,y)\in u\times u:x\neq
  y\}$. 
\item The Cantor space $\can$ of all infinite sequences with values 0 and 1
is equipped with the natural product topology and the group operation of
coordinate-wise addition $+$ modulo 2.   
\item Ordinal numbers will be denoted be the lower case initial letters of
  the Greek alphabet $\alpha,\beta,\gamma,\delta,\vare,\zeta$ as well as
  $\xi$. Finite ordinals (non-negative integers) will be denoted by letters
  $a,b,c,d,i,j,k,\ell,m,n,M$ and $\iota$. 
\item The Greek letters $\kappa,\lambda$ will stand for uncountable
  cardinals. 
\item For a forcing notion $\bbP$, all $\bbP$--names for objects in
  the extension via $\bbP$ will be denoted with a tilde below (e.g.,
  $\name{\tau}$, $\name{X}$), and $\name{G}_\bbP$ will stand for the
  canonical $\bbP$--name for the generic filter in $\bbP$.
\end{enumerate}
We fully utilize the algebraic properties of $(\can,+)$, in particular the fact 
that all elements of $\can$ are self-inverse.  

\section{The splitting rank}
In this section we remind some basic facts from \cite[Section 1]{Sh:522}
concerning a rank (on models with countable vocabulary) which will be used
in the construction of a forcing notion in the fourth section. This rank and
relevant proofs were also presented in \cite[Section 2]{RoSh:1138}.

Let $\lambda$ be a cardinal and $\bbM$ be a model with the universe
$\lambda$ and a countable vocabulary $\tau$.

\begin{definition}
\label{defofrank}
\begin{enumerate}
\item By induction on ordinals $\delta$, for finite non-empty sets
  $w\subseteq\lambda$ we define when $\rk(w,\bbM)\geq \delta$. Let
  $w=\{\alpha_0,\ldots,\alpha_n\} \subseteq\lambda$, $|w|=n+1$.  
 \begin{enumerate}
\item[(a)] $\rk(w)\geq 0$ if and only if for every quantifier free formula 
    $\varphi\in \cL(\tau)$ and each $k\leq n$, if 
$\bbM\models \varphi[\alpha_0,\ldots,\alpha_k,\ldots,\alpha_n]$ then  the set 
\[\big\{\alpha\in \lambda:\bbM\models \varphi[\alpha_0,\ldots,\alpha_{k-1},
  \alpha,\alpha_{k+1}, \ldots,\alpha_n]\big\}\]  is 
uncountable;  
\item[(b)] if $\delta$ is limit, then $\rk(w,\bbM)\geq\delta$ if and only if 
  $\rk(w,\bbM)\geq\gamma$ for all $\gamma<\delta$;  
\item[(c)] $\rk(w,\bbM)\geq\delta+1$ if and only if for every quantifier
  free  formula $\varphi\in \cL(\tau)$ and each $k\leq n$, if $\bbM\models  
  \varphi[\alpha_0,\ldots,\alpha_k,\ldots,\alpha_n]$ then there is
  $\alpha^*\in\lambda\setminus w$ such that 
\[\rk(w\cup\{\alpha^*\},\bbM)\geq \delta\quad\mbox{ and }\quad \bbM\models  
  \varphi[\alpha_0,\ldots,\alpha_{k-1},\alpha^*,\alpha_{k+1},\ldots,\alpha_n].\] 
 \end{enumerate}
\item The rank $\rk(w,\bbM)$ of a finite non-empty set 
$w\subseteq\lambda$ is defined by: 
 \begin{itemize}
\item $\rk(w,\bbM)=-1$ if $\neg (\rk(w,\bbM)\geq 0)$, and
\item $\rk(w,\bbM)=\infty$ if $\rk(w,\bbM)\geq \delta$ for all ordinals 
  $\delta$, and
\item for an ordinal $\delta$: $\rk(w,\bbM)=\delta$ if $\rk(w,\bbM)\geq 
  \delta$ but $\neg(\rk(w,\bbM)\geq\delta+1)$.
  \end{itemize}
\end{enumerate}
\end{definition}

\begin{definition}
For an ordinal $\vare$ and a cardinal $\lambda$ let ${\rm 
  NPr}_\vare(\lambda)$  be the following statement: ``there is a model 
$\bbM^*$ with the universe $\lambda$ and a countable vocabulary $\tau^*$
such that  $\sup\{\rk(w,\bbM^*):\emptyset\neq w\in [\lambda]^{<\omega}
\}<\vare$.''    

${\rm Pr}_\vare(\lambda)$ is the negation of ${\rm
  NPr}_\vare(\lambda)$.
\end{definition}

\begin{observation}
\label{obsonrk}
If $\lambda$ is uncountable and ${\rm NPr}_\vare(\lambda)$, then
  there is a model $\bbM^*$ with the universe $\lambda$ and a countable
  vocabulary $\tau^*$ such that 
\begin{itemize}
\item $\rk(\{\alpha\},\bbM^*)\geq 0$ for all $\alpha\in \lambda$ and  
\item  $\rk(w,\bbM^*)<\vare$ for every finite non-empty set $w\subseteq 
  \lambda$. 
\end{itemize}
\end{observation}

\begin{proposition}
[See {\cite[Claim 1.7]{Sh:522}} and/or {\cite[Proposition
  2.6]{RoSh:1138}}] 
\label{cl1.7-522}
\ \ 
\begin{enumerate}
\item ${\rm NPr}_1(\omega_1)$.
\item If ${\rm NPr}_\vare(\lambda)$, then ${\rm
    NPr}_{\vare+1}(\lambda^+)$. 
\item If ${\rm NPr}_\vare(\mu)$ for $\mu<\lambda$ and
  $\cf(\lambda)=\omega$, then ${\rm NPr}_{\vare+1}(\lambda)$. 
\end{enumerate}
\end{proposition}

\begin{proposition}
[See {\cite[Conclusion 1.8]{Sh:522}} and/or {\cite[Proposition 2.7]{RoSh:1138}}]
\label{522.1.8}  
Assume $\beta<\alpha<\omega_1$, $\bbM$ is a model with a countable
vocabulary $\tau$ and the universe $\mu$,  $m,n<\omega$, $n>0$, $A\subseteq 
\mu$ and $|A|\geq \beth_{\omega\cdot\alpha}$. Then there is $w\subseteq A$
with $|w|=n$ and $\rk(w,\bbM)\geq \omega\cdot \beta+m$ \footnote{``
    $\cdot$ '' stands for the ordinal multiplication}. 
\end{proposition}

\begin{definition}
\label{deflam}
Let $\lambda_{\omega_1}$ be the smallest cardinal $\lambda$ such that ${\rm   
  Pr}_{\omega_1}(\lambda)$. 
\end{definition}

\begin{corollary}
  \begin{enumerate}
\item If $\alpha<\omega_1$, then ${\rm NPr}_{\omega_1}(\aleph_\alpha)$.
\item ${\rm  Pr}_{\omega_1}(\beth_{\omega_1})$ holds true. 
\item $\aleph_{\omega_1}\leq\lambda_{\omega_1} \leq \beth_{\omega_1}$.      
  \end{enumerate}
\end{corollary}

\begin{corollary}
  [See {\cite[Proposition 2.10 and Corollary 2.11]{RoSh:1138}}]
 \label{lamCoh}
 Let $\mu=\beth_{\omega_1}\leq\kappa$. If $\bbP$ is a ccc forcing
 notion, then $\forces_{\bbP}{\rm Pr}_{\omega_1}(\mu)$. In particular,
 if $\bbC_\kappa$ be the forcing notion adding $\kappa$ Cohen reals,
 then $\forces_{\bbC_\kappa}\lambda_{\omega_1}\leq\mu\leq\con$.  
\end{corollary}

\section{Spectrum of translation non-disjointness}

We want to analyze sets with many non-disjoint translations in more detail,
restricting ourselves to $\Sigma^0_2$ subsets of $\can$. In this section we
will keep the following assumptions.

\begin{hypothesis}
\label{hyp1}
Let $\bar{T}=\langle T_n: n<\omega \rangle$, where each $T_n\subseteq 
  {}^{\omega>} 2$ is a tree with no maximal nodes (for $n<\omega$). Let 
 $B=\bigcup\limits_{n<\omega} \lim(T_n)$.
\end{hypothesis}

\begin{definition}
  \label{basedef}
  \begin{enumerate}
\item Let $\cL$ consist of all non-empty sets $u\subseteq \fstwo$ such that 
    $u\subseteq {}^\ell 2$ for some $\ell=\ell(u)<\omega$.
\item {\em A simple base\/} is a (strict) partial order $\cO=(\cO,\prec)$
  such that $\cO\subseteq \cL$ and for $u,u'\in\cO$:
  \begin{enumerate}
\item[(a)] if $u\prec u'$ then $\ell(u)<\ell(u')$ and $u= \{\eta\rest \ell(u): 
    \eta\in u'\}$, 
\item[(b)] there is a $v\in\cO$ such that $u\prec v$, 
\item[(c)] if $\rho\in {}^{\ell(u)}2$ then $u+\rho\in\cO$, and if $\rho\in
  {}^{\ell(u')} 2$ and $u\prec u'$  then $u+\rho\rest \ell(u)\prec u'+\rho$.
  \end{enumerate}
\item Let $(\cO,\prec)$ be a simple base. An {\em $\cO$--tower\/} is a 
  $\prec$--increasing sequence $\bar{u}= \langle
  u_n:n<\omega \rangle\subseteq \cO$ (so $u_n\prec u_{n+1}$ for all
  $n<\omega$). {\em The cover of an $\cO$--tower $\bar{u}$\/} is the set
  $\cC(\bar{u})\stackrel{\rm def}{=}\big\{\eta\in\can:\big(\forall
  n<\omega\big) \big(\eta\rest\ell(u_n)\in u_n\big)\big\}$.
\item  {\em An indexed base\/} is a sequence $\bar{\cO}=\langle \cO_i:i<i^*
  \rangle$ where $0<i^*\leq \omega$ and each $\cO_i$ is a simple base.
  \end{enumerate}
\end{definition}

\begin{definition}
\label{otps}
Let $\bar{\cO}=\langle \cO_i:i<i^* \rangle$ be an indexed base.
\begin{enumerate}
\item We say that two translations $B+x$ and $B+y$ of the
  set\footnote{Remember Assumptions \ref{hyp1}} $B$ (for $x,y\in \can$)
  have {\em $\bar{\cO}$--large intersection\/} if for some
  $\langle\bar{u}_i:i<i^*\rangle$ for every $i<i^*$ we have:
  \begin{itemize}
\item $\bar{u}_i$ is an $\cO_i$--tower,
\item for some $n_1,n_2<\omega$,
  \[\cC(\bar{u}_i)\subseteq \big(\lim(T_{n_1})+x\big)\cap 
    \big(\lim(T_{n_2})+y\big),\]     
\item $\cC(\bar{u}_i)\cap \cC(\bar{u}_j)=\emptyset$ whenever $j<i^*$, $j\neq i$.
\end{itemize}
In the above situation we may also say that {\em $(B+x)\cap (B+y)$ is
  $\bar{\cO}$--large.} 
\item We say that $B$ is {\em perfectly orthogonal to $\bar{\cO}$--small\/} (or
  a $\bar{\cO}$--{\bf pots}--set) if there is a perfect set $P\subseteq \can$
  such that the translations $B+x$, $B+y$ have a $\bar{\cO}$--large
  intersection for all $x,y\in P$. \\ 
The set $B$ is {\em an $\bar{\cO}$--{\bf npots}--set\/} if it is not
$\bar{\cO}$--{\bf pots}.   
\item We say that $B$ has {\em $\lambda$ many pairwise
    $\bar{\cO}$--nondisjoint translations\/} if for some set
  $X\subseteq\can$ of cardinality $\lambda$, for all $x,y\in X$ the
  translations $B+x$, $B+y$ have a $\bar{\cO}$--large intersection.
\item We define the {\em spectrum of translation
    $\bar{\cO}$--nondisjointness of $B$} as  
  \[\begin{array}{ll}
    \stnd_{\bar{\cO}} (B)=\big\{(x,y)\in\can\times\can: &\mbox{the translations
  }B+x,B+y\\
&\mbox{have a $\bar{\cO}$--large intersection. }\ \ \big\}.
\end{array} \]  
\end{enumerate} 
\end{definition}

\begin{example}
  \label{basex}
  \begin{enumerate}
  \item Let $3\leq \iota\leq\omega$. Put $\cO^0=\{u\in\cL: |u|=1\}$ and let
    a relation $\prec^0$ be defined by:
    
    $u\prec^0v$ if and only if $\ell(u)<\ell(v)\ \wedge\  u=\{\eta\rest
    \ell(u):\eta\in v\}$. 

    \noindent Then $(\cO^0,\prec^0)$ is a simple base and $\bar{\cO}^\iota=
    \langle \cO^0:i<\iota \rangle$ is an indexed base. Two translations
    $B+x$ and $B+y$ of the set $B$ (for $x,y\in \can$) have
    $\bar{\cO}^\iota$--large intersection if and only if $(B+x)\cap (B+y)$
    has at least $\iota$ members. 
  \item Let $\cO^\per=\{u\in \cL:|u|\geq 3\}$ and let a relation
    $\prec^\per$ be defined by
    
    $u\prec^\per v$ if and only if

    $ u=\{\eta\rest
      \ell(u):\eta\in v\}\  \wedge\ (\forall\nu\in u)(|\{\eta\in
      v:\nu\vtl\eta\}|\geq 2)$.\\
    Then $(\cO^\per,\prec^\per)$ is a simple base and $\bar{\cO}^\per=
    \langle \cO^\per\rangle$ is an indexed base. Two translations $B+x$ and
    $B+y$ of the set $B$ (for $x,y\in \can$) have   $\bar{\cO}^\per$--large
    intersection if and only if $(B+x)\cap (B+y)$ is uncountable.    
      \end{enumerate}
\end{example}

\begin{proposition}
  \label{proptostart}
  Let $\bar{\cO}$ be an indexed base and let $\bar{T},B$ be as in Assumptions
  \ref{hyp1}. 
\begin{enumerate}
\item The set $B$ is a $\bar{\cO}$--{\bf pots}--set if and only if  there
  is a perfect set $P\subseteq\can$ such that $P\times
  P\subseteq\stnd_{\bar{\cO}} (B)$.  
\item The set $\stnd_{\bar{\cO}}(B)$ is $\Sigma^1_1$.
\item Let $\con<\lambda\leq\mu$ and let $\bbC_\mu$ be the forcing notion
  adding $\mu$  Cohen reals. Then, remembering Definition \ref{otps}(2),  
\[\begin{array}{l}
\forces_{\bbC_\mu}\mbox{`` if $B$ has $\lambda$
    many pairwise $\bar{\cO}$--nondisjoint translations,}\\
\qquad\quad\mbox{then $B$ is a $\bar{\cO}$--{\bf pots}--set ''}. 
\end{array}\]
\item Assume ${\rm Pr}_{\omega_1}(\lambda)$. If $B$ has $\lambda$ many
  pairwise $\bar{\cO}$--nondisjoint translations, then it is an $\bar{\cO}$--{\bf
    pots}--set. 
\end{enumerate}
\end{proposition}

\begin{proof}
(1,2)\quad Straightforward; in evaluation of the complexity of
$\stnd_{\bar{\cO}}(B)$ note that for $\cO_i$--towers $\bar{u}_i=\langle 
u^i_n: n<\omega\rangle$, $x\in\can$ and $k<\omega$:

$\cC(\bar{u}_i)\subseteq \lim(T_k)+x$ if and only if $(\forall
  n<\omega)(u^i_n\subseteq T_k+x)$, and

  $\cC(\bar{u}_{i_1})  \cap \cC(\bar{u}_{i_2})=\emptyset$ if and only if
  $(\exists \ell<\omega)(\forall n_1, n_2>\ell)(u_{n_1}^{i_1}\rest \ell
  \cap u_{n_2}^{i_2}\rest \ell=\emptyset)$.   
\smallskip 

\noindent (3)\quad This is a consequence of (1,2) above and Shelah
\cite[Fact 1.16]{Sh:522}. 
\medskip

\noindent (4)\quad By \cite[Claim 1.12(1)]{Sh:522}. 
\end{proof}
\bigskip

To carry out our arguments we need to assume that our indexed base $\bar{\cO}$
satisfies some additional properties.

\begin{definition}
  \label{nicedef}
  An indexed base $\bar{\cO}=\langle\cO_i:i<i^*\rangle$ is {\em nice\/}
  if it satisfies the following demands (i)--(v). 
\begin{enumerate}
\item[(i)] Either $i^*\geq 6$ or for some $i<i^*$ we
  have 
\[\big(\forall u\in\cO_i\big)\big(\exists v\in\cO_i\big)(u\prec v\
  \wedge\ |v|\geq 6\big).\]
\item[(ii)] If $i<i^*$, $u\prec_i v\prec_i v'\prec_i v''$, and $\ell(v)\leq
  \ell\leq \ell(v')$, then $\{\eta\rest\ell:\eta\in v'\}\in\cO_i$ and
  $u\prec_i \{\eta\rest\ell:\eta\in v'\}\prec_i v''$.
\item[(iii)] If $i<i^*$, $u\prec_i v$, $\ell(v)<\ell$ and $v'\subseteq
  {}^\ell 2$ is such that for each $\nu\in v$ the set  $\{\eta\in
  v':\nu\vtl\eta\}$ has exactly one element, then $v'\in\cO_i$ and $u\prec_i
  v'$. 
\item[(iv)] Suppose $u\prec_i v$ and $u'\subseteq u$ is such that
  $u'\in\cO_i$. Let $v'=\{\eta\in v:\eta\rest\ell(u)\in u'\}$. Then
  $v'\in\cO_i$ and $u'\prec_i v'$. 
\item[(v)]  If $i^*=\omega$, then for each $i<i^*$ there are infinitely many
  $j<i^*$ such that $\cO_i=\cO_j$.
\end{enumerate}
\end{definition}

\begin{observation}
 The indexed bases $\bar{\cO}^\iota$ and $\bar{\cO}^\per$ introduced in
 Example \ref{basex} are nice.
\end{observation}

\begin{proposition}
  \label{niceprop}
Suppose an indexed base $\bar{\cO}=\langle\cO_i:i<i^*\rangle$ is
nice. Then:
  \begin{enumerate}
\item[$(\circledast)$]  If $2\leq K<\omega$ and $\bar{u}^k$ (for $k< 
    K$) is an  $\cO_{i(k)}$--tower for some $i(k)< i^*$, then there are 
    $\cO_{i(k)}$--towers $\bar{v}^k=\langle v^k_n:n<\omega\rangle$ (for 
    $k< K$) such that   
  \begin{itemize}
\item $\cC(\bar{v}^k)= \cC(\bar{u}^k)$, $v^k_0=u^k_0$ and 
\item $\bigcap\limits_{k\in K}\{\ell(v^k_n):n<\omega\}$ is infinite. 
  \end{itemize}
  \end{enumerate}
\end{proposition}

\begin{proof}
  Induction on $K$. For $K=2$ we proceed as follows. Let $\bar{u}^0$ be an
  $\cO_{i(0)}$--tower  and $\bar{u}^1$ be an $\cO_{i(1)}$--tower. Choose
  inductively a sequence $\langle n_k:k<\omega\rangle$ so that
  \begin{itemize}
  \item $5<n_0<n_1<n_2<n_3<\ldots$,
  \item $\ell(u^1_5)<\ell(u^0_{n_0})$,
  \item if $\ell(u^1_j)\leq \ell(u^0_{n_k})<\ell(n^1_{j+1})$, then
    $\ell(n^1_{j+5})\leq \ell(n^0_{k+1})$.
  \end{itemize}
  For $k<\omega$ let $j(k)$ be such that $\ell(u^1_{j(k)})\leq
  \ell(u^0_{n_k}) <\ell(n^1_{j(k)+1})$. Put $v_k=\{\eta\rest
  \ell(u^0_{n_k}): \eta\in u^1_{j(k)+1}\}$. By \ref{nicedef}(ii), $v_k\in
  \cO_{i(1)}$ and $u^1_{j(k)-1}\prec_{i(1)} v_k\prec_{i(1)}
  u^1_{j(k)+2}$. The rest should be clear.  
\end{proof}

For the rest of this section we will be assuming the following. 

\begin{hypothesis}
\label{hyp2}
  \begin{enumerate}
\item $\bar{T}=\langle T_n: n<\omega \rangle$, $B$ are as in Assumptions
  \ref{hyp1},
\item $\bar{\cO}=\langle\cO_i:i<i^*\rangle$ is a nice indexed base. Also,
  $\cO_i=(\cO_i,\prec_i)$,
\item there are distinct $x,y\in\can$ such that $(B+x)\cap (B+y)$ is
  $\bar{\cO}$--large. 
  \end{enumerate}
\end{hypothesis}

\begin{definition}
  \label{mtkDef}
Let $\Mtk$ consist of all tuples 
\[\bm=(\ell^\bm,\iota^\bm,u^\bm,\bar{h}^\bm,\bar{g}^\bm)=(\ell,\iota,u, \bar{h},
  \bar{g})\]   
such that:
\begin{enumerate}
\item[(a)] $0<\ell<\omega$, $u\subseteq {}^\ell 2$ and $2\leq |u|$, and
  $\iota=i^*$ if $i^*<\omega$, and $3\leq \iota<\omega$ otherwise;  
\item[(b)] $\bar{g}=\langle g_i:i<\iota\rangle$, where\footnote{remember
    $u^{\langle 2\rangle}=\{(\eta,\nu)\in u\times u: \eta\neq \nu\}$}
  $g_i:u^{\langle 2\rangle}\longrightarrow \cO_i$ is such that
  $g_i(\eta,\nu)= g_i(\nu,\eta)$ and $\ell\big(g_i(\eta,\nu)\big)=\ell$ for
  each  $(\eta,\nu)\in u^{\langle 2\rangle}$;  
\item[(c)] if $(\eta,\nu)\in u^{\langle 2\rangle}$ and $i<i'<\iota$,
  then $g_i(\eta,\nu)\cap g_{i'}(\eta,\nu)=\emptyset$, 
\item[(d)] $\bar{h}=\langle h_i:i<\iota\rangle$, where $h_i:u^{\langle  
    2\rangle}\longrightarrow \omega$; 
\item[(e)] for each $(\eta,\nu)\in u^{\langle 2\rangle}$, if $\sigma\in
  g_i(\eta,\nu)$ then $\eta+\sigma\in T_{h_i(\eta,\nu)}$.
\end{enumerate}
\end{definition}

\begin{definition}
\label{traDef}
Assume $\bm=(\ell,\iota,u,\bar{h},\bar{g})\in \Mtk$ and $\rho\in {}^\ell
2$. We define $\bm+\rho=(\ell',\iota',u',\bar{h}',\bar{g}')$ by
\begin{itemize}
\item $\ell'=\ell$, $\iota'=\iota$, $u'=\{\eta+\rho:\eta\in u\}$,
\item $\bar{g}'=\langle g'_i:i<\iota'\rangle$, where $g'_i:(u')^{\langle
    2\rangle} \longrightarrow \cO_i:(\eta+\rho,\nu+\rho)\mapsto
  g_i(\eta,\nu)+\rho$, 
\item $\bar{h}'=\langle h'_i:i<\iota'\rangle$, where $h_i':(u')^{\langle
    2\rangle} \longrightarrow\omega$ are such that $h'_i(\eta+\rho,\nu+\rho)
  = h_i(\eta,\nu)$ for $(\eta,\nu)\in u^{\langle 2\rangle}$.
\end{itemize}
Also if $\rho\in\can$, then we set $\bm+\rho=\bm+(\rho\rest\ell)$. 
\end{definition}

\begin{observation}
  \begin{enumerate}
\item If $\bm\in \Mtk$ and $\rho\in {}^{\ell^\bm}2$, then $\bm+\rho
  \in\Mtk$. 
\item For each $\rho\in\can$ the mapping $\Mtk\longrightarrow \Mtk:\bm
  \mapsto\bm+\rho$ is a bijection.
  \end{enumerate}
\end{observation}

\begin{definition}
\label{extDef}
Assume $\bm,\bn\in\Mtk$. We say that {\em $\bn$ strictly extends $\bm$\/}
($\bm\sqsubset \bn$ in short) if and only if:
\begin{itemize}
\item $\ell^\bm< \ell^\bn$,  $\iota^\bm\leq \iota^\bn$, 
  $u^\bm=\{\eta\rest\ell^\bm:\eta\in u^\bn\}$, and
\item for every $(\eta,\nu)\in (u^\bn)^{\langle 2\rangle }$ such that
  $\eta\rest \ell^\bm \neq \nu\rest\ell^\bm$ and each $i<\iota^\bm$ we have   
  \begin{itemize}
  \item  $g^\bm_i(\eta\rest \ell^\bm,\nu\rest \ell^\bm)\prec g^\bn_i(\eta,
    \nu)$, and
  \item  $h^\bm_i(\eta\rest \ell^\bm,\nu\rest \ell^\bm)= h^\bn_i(\eta,\nu)$.   
  \end{itemize}
\end{itemize}
\end{definition}

\begin{definition}
  \label{ndrkdef}
\begin{enumerate}
\item By induction on ordinals $\alpha$ we define
$D^{\bar{T}}(\alpha)\subseteq\Mtk$. We declare that: 
\begin{itemize}
\item $D^{\bar{T}}(0)=\Mtk$, 
\item if $\alpha$ is a limit ordinal, then $D^{\bar{T}}(\alpha)=
  \bigcap\limits_{\beta<\alpha} D^{\bar{T}}(\beta)$, 
\item if $\alpha=\beta+1$, then $D^{\bar{T}}(\alpha)$ consists of all
  $\bm\in\Mtk$ such that for each for each $\nu\in u^\bm$ there is an $\bn\in
  \Mtk$ satisfying
  \begin{itemize}
  \item $\bm\sqsubset\bn$ and $\bn\in D^{\bar{T}}(\beta)$, and if
    $i^*=\omega$ then $\iota^\bm<\iota^\bn$, and 
  \item the set $\{\eta\in u^\bn:\nu\vtl \eta\}$ has at least two elements
  \end{itemize}
\end{itemize}
\item We define a function\footnote{$\ndrk$ stands for\quad  {\bf n}on{\bf
      d}isjointness {\bf r}an{\bf k}} $\ndrk^{\bar{T}}_{\bar{\cO}}=\ndrk:\Mtk
  \longrightarrow {\rm ON}\cup\{\infty\}$ as follows.\\
    If $\bm\in D^{\bar{T}}(\alpha)$ for all ordinals $\alpha$,  then 
    we say that $\ndrk(\bm)=\infty$.\\ 
    Otherwise, $\ndrk(\bm)$ is the first ordinal $\alpha$ for which 
    $\bm\notin D^{\bar{T}}(\alpha+1)$. 
\item We also define 
\[\NDRK_{\bar{\cO}}(\bar{T})=\NDRK(\bar{T})=\sup\{\ndrk(\bm)+1:\bm\in\Mtk\}.\]  
\end{enumerate}
\end{definition}

\begin{lemma}
\label{lemonrk}
  \begin{enumerate}
  \item The relation $\sqsubset$ is a strict partial order on $\Mtk$.
  \item If $\bm,\bn\in\Mtk$ and $\bm\sqsubset\bn$ and $\bn\in
    D^{\bar{T}}(\alpha)$, then $\bm\in D^{\bar{T}}(\alpha)$.
\item If $\alpha<\beta$ then  $D^{\bar{T}}(\beta)\subseteq
  D^{\bar{T}}(\alpha)$. Hence for $\bm\in\Mtk$, $\bm\in D^{\bar{T}}(\alpha)$
  if and only if $\alpha\leq\ndrk(\bm)$. 
\item If $\bm\in\Mtk$ and $\rho\in\can$ then $\ndrk(\bm)=\ndrk(\bm+\rho)$. 
\item If $\bm\in\Mtk$ and $\ndrk(\bm)\geq \omega_1$, then
  there is an $\bn\in\Mtk$ such that $\bm\sqsubset\bn$, $|\{\eta\in
  u^\bn:\nu\vtl\eta\}|\geq 2$ for each $\nu\in u^\bm$, if $i^*=\omega$ then
  $\iota^\bm<\iota^\bn$, and $\ndrk(\bn)\geq\omega_1$.
 \item If $\bm\in\Mtk$ and $\infty>\ndrk(\bm)=\beta>\alpha$, then there is
  $\bn\in \Mtk$ such that $\bm\sqsubset\bn$ and $\ndrk(\bn)=\alpha$. 
\item If $\NDRK(\bar{T})\geq \omega_1$, then $\NDRK(\bar{T})=\infty$. 
\item Assume $\bm\in\Mtk$ and $u'\subseteq u^\bm$, $|u'|\geq 2$. Put
$\ell'=\ell^\bm$, $\iota'=\iota^\bm$, and for $i<\iota'$ let
$h'_i=h^\bm_i\rest (u')^{\langle 2\rangle}$ and $g'_i=g^\bm_i \rest
(u')^{\langle 2\rangle}$.  Let  $\bm\rest
u'=(\ell',u',i',\bar{h}',\bar{g}')$. Then $\bm \rest u'\in \Mtk$  and
$\ndrk(\bm)\leq \ndrk(\bm\rest u')$.   
  \end{enumerate}
\end{lemma}

\begin{proof}
Exactly the same as for \cite[Lemma 3.10]{RoSh:1138}.
\end{proof}

\begin{proposition}
\label{eqnd}
For a nice indexed base $\bar{\cO}$ the following conditions (a) -- (d)
are equivalent.  
\begin{enumerate}
\item[(a)] $\NDRK_{\bar{\cO}}(\bar{T})\geq \omega_1$.
\item[(b)] $\NDRK_{\bar{\cO}}(\bar{T})=\infty$.
\item[(c)] $B$ is perfectly orthogonal to $\bar{\cO}$--small (see \ref{otps}(2)).
\item[(d)] In some ccc forcing extension, the set $B$ has
  $\lambda_{\omega_1}$ many pairwise $\bar{\cO}$--nondisjoint translations
  (see \ref{otps}(3)). 
\end{enumerate}
\end{proposition}

\begin{proof}
The proof follows closely the lines of \cite[Proposition 3.11]{RoSh:1138}. 

 \noindent ${\rm (c)} \Rightarrow {\rm (d)}$\quad Assume (c) and let
 $P\subseteq \can$ be a perfect set such that the translations $B+x$, $B+y$
 have $\bar{\cO}$--large intersection for all $x,y\in P$. Let
 $\kappa=\beth_{\omega_1}$. By Corollary \ref{lamCoh},
 $\forces_{\bbC_\kappa} \lambda_{\omega_1}\leq \con$. By Proposition  
\ref{proptostart}(2), the formula ``$P\times P\subseteq \stnd_{\bar{\cO}}
(B)$'' is $\Pi^1_2$, so it holds in the forcing extension by
$\bbC_\kappa$. Now we easily conclude (d).  
\medskip

\noindent ${\rm (d)} \Rightarrow {\rm (a)}$\quad Assume (d) and let $\bbP$
be the ccc forcing notion witnessing this assumption, $G\subseteq\bbP$ be
generic over $\bV$. Let us work in $\bV[G]$. 

Let $\langle \eta_\alpha:\alpha<\lambda_{\omega_1}\rangle$ be
a sequence of distinct elements of $\can$ such that
\[\big(\forall \alpha<\beta<\lambda_{\omega_1}\big)\big((B+\eta_\alpha)\cap 
(B+\eta_\beta)\mbox{ is $\bar{\cO}$--large }\big).\]
Remember Definition \ref{basedef}(3): an $\cO_i$--tower is an 
$\prec_i$--increasing sequence $\bar{u}=\langle u_n:n<\omega\rangle$ and its 
cover $\cC(\bar{u})$ is the set $\{\eta\in\can: (\forall
n<\omega)(\eta\rest\ell(u_n)\in  u_n)\}$.

Let $\tau=\{R_\bm:\bm\in\Mtk\}$ be a vocabulary where each $R_\bm$ is a
$|u^\bm|$--ary relational symbol. Let $\bbM=\big(\lambda_{\omega_1},
\big\{R^\bbM_\bm\big\}_{\bm\in \Mtk}\big)$ be the model in the vocabulary
$\tau$, where for $\bm=(\ell,\iota,u,h,g)\in\Mtk$ the relation $R_\bm^\bbM$ is
defined by
\[\begin{array}{ll}
R^\bbM_\bm=&\Big\{(\alpha_0,\ldots,\alpha_{|u|-1})\in
(\lambda_{\omega_1})^{|u|}:\{\eta_{\alpha_0}\rest \ell,\ldots,
             \eta_{|u|-1}\rest \ell\} =u \mbox{ and}\\ 
&\quad \mbox{for each distinct }j_1,j_2<|u|\mbox{ and every } i<\iota\\
&\quad \mbox{there is an $\cO_i$--tower }\bar{u}^i(j_1,j_2)=\langle
 u_n^i(j_1,j_2)\!: n<\omega\rangle\mbox{ such that }\ \\             
&\quad g_i(\eta_{\alpha_{j_1}}\rest\ell, \eta_{\alpha_{j_2}}\rest 
\ell)=u_0^i(j_1,j_2)\mbox{ and }\cC\big(\bar{u}^i(j_1,j_2)\big) \mbox{ is
     included in }\\
&\quad [\lim(T_{h_i(\eta_{\alpha_{j_1}}\rest \ell,\eta_{\alpha_{j_2}}\rest\ell)})+ 
\eta_{\alpha_{j_1}}] \cap [\lim(T_{h_i(\eta_{\alpha_{j_2}},
            \eta_{\alpha_{j_1}})}) +\eta_{\alpha_{j_2}}] \Big\}.
\end{array}\]

\begin{claim}
\label{cl5}
\begin{enumerate}
\item If $\alpha_0,\alpha_1,\ldots,\alpha_{j-1}<\lambda_{\omega_1}$ are
  distinct, $j\geq 2$, then for infinitely many $k<\omega$ there is
  $\bm\in \Mtk$ such that  
\[\ell^\bm=k,\quad u^\bm=\{\eta_{\alpha_0}\rest k, \ldots,
\eta_{\alpha_{j-1}}\rest k\}\quad \mbox{ and }\quad \bbM\models
R_\bm[\alpha_0,\ldots,\alpha_{j-1}].\]

\item Assume that  $\bm\in\Mtk$, $j<|u^\bm|$, 
  $\alpha_0,\alpha_1,\ldots,\alpha_{|u^{\bm}|-1} < \lambda_{\omega_1}$ and
  $\alpha^*< \lambda_{\omega_1}$ are all pairwise distinct and such that
  \begin{itemize}
  \item   $\bbM\models R_{\bm}[\alpha_0,\ldots,\alpha_j, \ldots,
    \alpha_{|u^{\bm}|-1}]$  and
  \item $\bbM\models R_{\bm}[\alpha_0,\ldots,\alpha_{j-1},\alpha^*,
    \alpha_{j+1}, \ldots \alpha_{|u^{\bm}|-1}]$.
  \end{itemize}
Then for infinitely many 
$k<\omega$ there is an $\bn\in \Mtk$ such that $\bm\sqsubset \bn$ 
and $\ell^\bn=k$, $u^\bn=\{\eta_{\alpha_0}\rest k, \ldots,
\eta_{\alpha_{|u^\bm|-1}}\rest k,\eta_{\alpha^*}\rest k\}$ and $\bbM\models
R_\bn[\alpha_0,\ldots,\alpha_{|u^\bm|-1},\alpha^*]$, and if $i^*=\omega$
then also $\iota^\bm<\iota^\bn$. 
\item If $\bm\in\Mtk$ and $\bbM\models R_\bm[\alpha_0,\ldots,
  \alpha_{|u^\bm|-1}]$, then  
\[\rk(\{\alpha_0,\ldots,\alpha_{|u^\bm|-1}\},\bbM)\leq
  \ndrk^{\bar{T}}_{\bar{\cO}}(\bm).\]   
\end{enumerate}
\end{claim}

\begin{proof}[Proof of the Claim]
  (1)\ It is a simpler version of the proof below.
\medskip

\noindent (2)\ By the definition of $R^\bbM_\bm$, since $\bbM\models
R_{\bm}[\alpha_0,\ldots,\alpha_{j-1},\alpha^*,\alpha_{j+1}, \ldots
\alpha_{|u^{\bm}|-1}]$ and $\bbM\models
R_{\bm}[\alpha_0,\ldots,\alpha_j, \ldots, \alpha_{|u^{\bm}|-1}]$, 
we may choose a sequence
\[\langle \bar{u}^i(j_1,j_2): (j_1,j_2)\in \big(|u^\bm|+1
\big)^{\langle 2\rangle}\ \wedge\ i<\iota^\bm\rangle\]
satisfying the following demands. Letting $\alpha_{|u^\bm|}=\alpha^*$, for
$(j_1,j_2)\in \big(|u^\bm|+1\big)^{\langle 2\rangle}$ and $i<\iota^\bm$:
\begin{itemize}
\item $\bar{u}^i(j_1,j_2)=\bar{u}^i(j_2,j_1)$ is a $\cO_i$--tower,
\item if $\{j_1,j_2\}\neq \{j,|u^\bm|\}$, then $u^i_0(j_1,j_2)= g_i^\bm(
  \eta_{\alpha_{j_1}}\rest \ell^\bm,  \eta_{\alpha_{j_2}}\rest \ell^\bm)$,
\item if $i_1<i_2<\iota^\bm$, then $\cC\big(\bar{u}^{i_1}(j_1,j_2)\big) \cap
  \cC\big(\bar{u}^{i_2}(j_1,j_2)\big)=\emptyset$, 
\item  if $\{j_1,j_2\}\neq \{j,|u^\bm|\}$, then
  $\cC\big(\bar{u}^i(j_1,j_2)\big)$ is included in 
\[ [\lim(T_{h^\bm_i(\eta_{\alpha_{j_1}}\rest \ell^\bm,
    \eta_{\alpha_{j_2}}\rest \ell^\bm)})+ \eta_{\alpha_{j_1}}] \cap 
  [\lim(T_{h_i^\bm(\eta_{\alpha_{j_2}}\rest \ell^\bm,
    \eta_{\alpha_{j_1}}\rest \ell^\bm)} ) +\eta_{\alpha_{j_2}}],\]
\item for some $N_i',N_i''$ we have 
\[\cC\big(\bar{u}^i(j,|u^\bm|)\big)\subseteq
  [\lim(T_{N_i'})+ \eta_{\alpha_j}] \cap [\lim(T_{N_i''})
  +\eta_{\alpha^*}].\] 
\end{itemize}
Since $\bar{\cO}$ is nice (and $\iota^\bm$ and $u^\bm$ are finite), we
may use \ref{niceprop}$(\circledast)$ and modify $\bar{u}^i(j_1,j_2)$
(without changing $u^i_0(j_1,j_2)$) and demand that the set 
\[A=\bigcap_{i<\iota^\bm} \bigcap_{j_1<j_2\leq |u^\bm|}\big\{\ell\big(
  u^i_n(j_1, j_2)\big):n\in\omega\big\}\] 
is infinite. Let $\ell_0\in A\setminus (\ell^\bm+1)$ be bigger than the
second element of $A\setminus (\ell^\bm+1)$ and such that
$\eta_{\alpha_{|u^\bm|}}\rest\ell_0 \neq\eta_{\alpha_j} \rest \ell_0$, and
$x\rest \ell_0\neq y\rest\ell_0$ whenever $x\in\cC(\bar{u}^{i_1}(j_1,j_2))$,
$y\in\cC(\bar{u}^{i_2}(j_1,j_2))$, $(j_1,j_2)\in \big(|u^\bm|+1
\big)^{\langle 2\rangle}$ and $i_1<i_2<\iota^\bm$.    

Let $\iota=\iota^\bm=i^*$ if $i^*<\omega$ and let $\iota=\iota^\bm+1$
otherwise. In the latter case we also have to choose
$\cO_{\iota^\bm}$--towers $\bar{u}^{\iota^\bm}(j_1,j_2)$, but to ensure the
demand \ref{mtkDef}(c) we will have to modify the already chosen towers
$\bar{u}^i(j_1,j_2)$ (for $i<\iota^\bm$). Fix $(j_1,j_2)\in \big(|u^\bm|+
1\big)^{\langle 2\rangle}$ for a moment.  Let 
\[K=\sum\big\{\big|u^i_n(j_1,j_2)\big|:\ell(u^i_n(j_1,j_2))=\ell_0\ \wedge\
  i<\iota^\bm\ \wedge\ n<\omega\big\}.\]
By \ref{nicedef}(v) and the assumptions on
$\langle\eta_\alpha:\alpha<\lambda_{\omega_1}\rangle$, there 
are infinitely many $\cO_{\iota^\bm}$--towers $\bar{v}^k$ such that their
covers are pairwise disjoint and included in
$\big(\lim(T_{k_1})+\eta_{\alpha_{j_1}}\big) \cap
\big(\lim(T_{k_2})+\eta_{\alpha_{j_2}}\big)$ for some $k_1,k_2$.

Choose $\ell(j_1,j_2)\in A\setminus (\ell_0+1)$ so large, that there are
more than $K+1$ many $k$'s for which the sets $\{\eta\rest\ell(j_1,j_2):
\eta \in v^k_n\}$ are pairwise disjoint (for large $n$) and
$\ell(v^k_5)<\ell(j_1,j_2)$ for all those $k$'s. For $i<i^\bm$ let
$n(i),m(i)$ be such that $\ell(u^i_{n(i)}(j_1,j_2))=\ell_0$ and
$\ell(u^i_{m(i)}(j_1,j_2))=\ell(j_1,j_2)$, and let $v^i\subseteq
u^i_{m(i)}(j_1,j_2)$ be such that for each $\nu\in u^i_{n(i)}(j_1,j_2)$ the
set $\{\eta\in v^i:\nu\vtl \eta\}$ has exactly one element. By
\ref{nicedef}(iii) we have
\[v^i\in\cO_i\quad\mbox{ and }\quad u^i_{n(i)-1}(j_1,j_2)\prec_i v^i.\] 
Using repeatedly \ref{nicedef}(iv) we may modify the towers
$\bar{u}^i(j_1,j_2)$ (for $i<i^\bm$) and demand that 
\begin{itemize}
\item for each $i<\iota^\bm$, for some $n^*(i)$,
  \[\ell(u^i_{n^*(i)}(j_1,j_2))=\ell(j_1,j_2)\ \mbox{ and }\
    \big|u^i_{n^*(i)}(j_1,j_2) \big|=\big|\{\eta\rest \ell_0: \eta\in
    u^i_{n^*(i)}(j_1,j_2)\}\big|.\]
\end{itemize}
Looking back at the towers $\bar{v}^k$, we may choose one,
$\bar{v}^{k^*}=\bar{v}(j_1,j_2)$, which has the property that for all large
$n$ 
\[\{\eta\rest \ell(j_1,j_2): \eta\in v_n(j_1,j_2)\}\cap
  \bigcup\big\{u^i_{n^*(i)}(j_1,j_2):i<\iota^\bm\} =\emptyset.\]

Now unfix $(j_1,j_2)$ and set $\ell=\max\{\ell(j_1,j_2):(j_1,j_2)\in
(|u^\bm|+1)^{\langle 2\rangle}\}$.

Suppose $j_1<j_2\leq |u^\bm|$ and let $n$ be such
that $\ell(v_{n-1}(j_1,j_2))<\ell\leq \ell(v_n(j_1,j_2))$. By
\ref{nicedef}(ii), we may let
\begin{itemize}
\item $u^{\iota^\bm}_0(j_1,j_2)=u^{\iota^\bm}_0(j_2,j_1)=\{\eta\rest \ell:\eta\in 
  v_n(j_1,j_2)\}$,
\item $u^{\iota^\bm}_m(j_1,j_2)=u^{\iota^\bm}_m(j_2,j_1)=
  v_{n+m}(j_1,j_2)$ for $m>0$,
\end{itemize}
getting a $\cO_{\iota^\bm}$--tower $\bar{u}^{\iota^\bm}(j_1,j_2)$. We also
fix $k(j_1,j_2), k(j_2,j_1)$ such that
\[\cC(\bar{u}^{\iota^\bm}(j_1,j_2))\subseteq
  \big(\lim(T_{k(j_1,j_2)}+\eta_{\alpha_{j_1}}\big) \cap
  \big(\lim(T_{k(j_2,j_1)} +\eta_{\alpha_{j_2}}\big).\]
If $i^*=\iota^\bm<\omega$, then the procedure leading to the choice of
$\bar{u}^{\iota^\bm}(j_1,j_2)$ is not present and we just let
$\ell=\min(A\setminus (\ell_0+1))$.
\medskip

Let $u=\big\{\eta_{\alpha_0}\rest \ell,\ldots, \eta_{\alpha_{|u^\bm|-1}}\rest
\ell, \eta_{\alpha^*}\rest\ell\big\}$. 

For each $i<\iota$ and $(j_1,j_2)\in \big(|u^\bm|+1)^{\langle 2\rangle}$ put  
$g_i(\eta_{\alpha_{j_1}}\rest\ell,\eta_{\alpha_{j_2}}\rest\ell)=
u^i_n(j_1,j_2)$, where $n$ is such that  $\ell\big(u^i_n(j_1,j_2)\big)=
\ell$. This defines $g_i:u^{\langle 2\rangle} \longrightarrow \cO_i$ for
$i<\iota$.  For $(\nu_1,\nu_2)\in u^{\langle 2\rangle}$ we also set 
\[h_i(\nu_1,\nu_2)=\left\{
  \begin{array}{ll}
h^\bm_i(\nu_1\rest\ell^\bm, \nu_2\rest\ell^\bm)
    &\mbox{ if }\nu_1\rest \ell^\bm\neq\nu_2\rest\ell^\bm,\ i<\iota^\bm, \\   
N_i' &\mbox{ if }\nu_1\vtl \eta_{\alpha_j},\  \nu_2\vtl\eta_{\alpha^*},\
       i<\iota^\bm, \\
N_i'' &\mbox{ if }\nu_1\vtl \eta_{\alpha^*},\  \nu_2\vtl\eta_{\alpha_j},\
       i<\iota^\bm, \\
k(j_1,j_2) &\mbox{ if }\nu_1\vtl \eta_{\alpha_{j_1}},\ \nu_2\vtl
      \eta_{\alpha_{j_2}},\ i=\iota^\bm<\iota.
  \end{array}\right. \]
It should be clear that $\bn=(\ell,\iota,u,g,h)\in\Mtk$ is as required.  
\medskip

\noindent (3)\quad By induction on $\beta$ we show that
\begin{quotation}
{\em for every\/}
$\bm\in\Mtk$ and {\em all\/} $\alpha_0,\ldots,\alpha_{|u^\bm|-1}<
\lambda_{\omega_1}$ such that $\bbM\models R_\bm[\alpha_0,\ldots,
\alpha_{|u^\bm|-1}]$: 

$\beta\leq \rk(\{\alpha_0,\ldots,\alpha_{|u^\bm|-1}\}, \bbM)$ implies
$\beta\leq \ndrk(\bm)$. 
\end{quotation}
\smallskip

\noindent {\sc Steps $\beta=0$ and $\beta$ is limit:}\quad Straightforward.  
\smallskip

\noindent {\sc Step $\beta=\gamma+1$:}\quad Suppose $\bm\in\Mtk$ and 
$\alpha_0,\ldots,\alpha_{|u^\bm|-1}<\lambda_{\omega_1}$ are such that
$\bbM\models R_\bm[\alpha_0,\ldots, \alpha_{|u^\bm|-1}]$ and
$\gamma+1\leq \rk(\{\alpha_0,\ldots,\alpha_{|u^\bm|-1}\}, \bbM)$. Let
$\nu\in u^\bm$, so $\nu=\eta_{\alpha_j}\rest \ell^\bm$ for some $j<|u^\bm|$.
Since $\gamma+1\leq \rk(\{\alpha_0,\ldots,\alpha_{|u^\bm|-1}\}, \bbM)$ we may
find $\alpha^*\in \lambda_{\omega_1}\setminus\{\alpha_0,\ldots,
\alpha_{|u^\bm|-1}\}$ such that
\[\bbM\models R_\bm[\alpha_0,\ldots,
  \alpha_{j-1},\alpha^*, \alpha_{j+1}, \ldots,\alpha_{|u^\bm|-1}]\]
and $\rk(\{\alpha_0,\ldots,\alpha_{|u^\bm|-1}, \alpha^*\}, \bbM)\geq
\gamma$. By clause (2) we may find $\bn\in\Mtk$ such that $\bm\sqsubset \bn$
and $u^\bn=\{\eta_{\alpha_0}\rest \ell^\bn,\ldots,\eta_{\alpha_{|u^\bm|-1}} 
\rest \ell^\bn, \eta_{\alpha^*}\rest\ell^\bn\}$, and if $i^*=\omega$ then
$\iota^\bm< \iota^\bn$, and $\bbM\models
R_\bn[\alpha_0,\ldots,\alpha_{|u^\bm|-1},\alpha^*]$. Then also 
$|\{\eta\in u^\bn:\nu\vtl\eta\}|\geq 2$. By the inductive hypothesis we have
also $\gamma\leq \ndrk(\bn)$.  Now we may easily conclude that $\gamma+1\leq  
\ndrk(\bm)$. 
\end{proof}

By the definition of $\lambda_{\omega_1}$, 
\begin{enumerate}
\item[$(\odot)$] $\sup\{\rk(w,\bbM):\emptyset\neq w\in
  [\lambda_{\omega_1}]^{<\omega} \}\geq\omega_1$ 
\end{enumerate}
Now, suppose that $\beta<\omega_1$. By $(\odot)$, there are distinct 
$\alpha_0,\ldots,\alpha_{j-1}<\lambda_{\omega_1}$, $j\geq 2$, such that
$\rk(\{\alpha_0,\ldots, \alpha_{j-1}\},\bbM)\geq \beta$. By Claim \ref{cl5}(1)
we may find $\bm\in\Mtk$ such that $\bbM\models
R_\bm[\alpha_0,\ldots,\alpha_{j-1}]$. Then by Claim \ref{cl5}(3) we also
have $\ndrk^{\bar{T}}_{\bar{\cO}}(\bm)\geq \beta$. Consequently,
$\NDRK(\bar{T})\geq \omega_1$.   

All the considerations above where carried out in $\bV[G]$. However, the
rank function $\ndrk^{\bar{T}}_{\bar{\cO}}$ is absolute, so we may also
claim that in $\bV$ we have $\NDRK_{\bar{\cO}}(\bar{T})\geq \omega_1$.   
\end{proof}

\section{The main result}
In this section we construct a forcing notion adding a sequence $\bar{T}$ of
subtrees of ${}^{\omega>} 2$ such that $\NDRK_{\bar{\cO}^6}(\bar{T})
<\omega_1$ and yet with many $\bar{\cO}$--nondisjoint translations (for a
nice $\bar{\cO}$). The sequence $\bar{T}$ will be added by finite
approximations, so we will need a finite version of Definition \ref{mtkDef}.  

\begin{definition}
  \label{fmtkDef}
Assume that 
\begin{enumerate}
\item[(a)] $0<n,M<\omega$, $\bar{t}=\langle t_m:m<M\rangle$, and  each $t_m$ 
  is a subtree of ${}^{n\geq} 2$ in which all terminal  branches are of
  length $n$, 
\item[(b)] $T_j\subseteq {}^{\omega>} 2$ (for $j<\omega$) are trees with no 
  maximal nodes, $\bar{T}=\langle T_j:j<\omega\rangle$ and $t_m=T_m\cap 
  {}^{n\geq} 2$ for $m<M$, 
\item[(c)] $\bM_{\bar{T},\bar{\cO}^6}$ is defined as in Definition 
  \ref{mtkDef} for $\bar{\cO}^6$ introduced in Example \ref{basex}(1).   
\end{enumerate}
We let $\fMtk$ consist of all tuples $\bm=(\ell^\bm,6, u^\bm,
\bar{h}^\bm, \bar{g}^\bm)\in \bM_{\bar{T},\bar{\cO}^6}$ such that
$\ell^\bm\leq n$ and $\rng(h^\bm_i)\subseteq M$ for each $i<6$.

The extension relation $\sqsubset$ on $\fMtk$ is inherited from
$\bM_{\bar{T},\bar{\cO}^6}$ (see Definition \ref{extDef}). 
\end{definition}

\begin{observation}
  \begin{enumerate}
\item The Definition of $\fMtk$ does not depend on the choice of $\bar{T}$,
  as long as the clause \ref{fmtkDef}(c) is satisfied.
  \item If $\bm\in \fMtk$ and $\rho\in {}^{\ell^\bm}2$, then $\bm+\rho \in\fMtk$
(remember Definition \ref{traDef}).  
  \end{enumerate}
\end{observation}

\begin{lemma}
[See {\cite[Lemma 2.3]{RoRy18}}.]
 \label{litlem}
Let $0<\ell<\omega$ and let $\cB\subseteq {}^\ell 2$ be a linearly
independent set of vectors (in $({}^\ell2,+)$ over $\bbZ_2$).
If $\cA\subseteq {}^\ell 2$, $|\cA|\geq 5$ and $\cA+\cA\subseteq \cB+\cB$, 
then for a unique $x\in {}^\ell 2$ we have $\cA+x\subseteq \cB$.  
\end{lemma}

\begin{theorem}
\label{522fortranslate}
Assume that an uncountable cardinal $\lambda$ satisfies ${\rm
  NPr}_{\omega_1}(\lambda)$ and suppose that $\bar{\cO}=\langle
\cO_i:i<i^*\rangle$ is a nice indexed base. Then there is a ccc forcing 
notion $\bbP$ of size $\lambda$ such that   
\[\begin{array}{l}
\forces_{\bbP}\mbox{`` for some $\Sigma^0_2$ $\bar{\cO}^6$--{\bf npots}--set  
    }B=\bigcup\limits_{n<\omega}\lim(T_n)\subseteq\can\mbox{ there is }\\
    \qquad\mbox{ a sequence } \langle\eta_\alpha:\alpha<\lambda\rangle
    \mbox{ of distinct elements of $\can$ such that}\\ 
\qquad \mbox{ all intersections }(\eta_\alpha+B)\cap
    (\eta_\beta+B)\mbox{ are $\bar{\cO}$--large for }\alpha,\beta<
    \lambda\mbox{ ''.} 
\end{array}\]
\end{theorem}

\begin{proof}
Fix a countable vocabulary $\tau=\{R_{n,\zeta}:n,\zeta< \omega\}$,  where
$R_{n,\zeta}$ is an $n$--ary relational symbol (for $n,\zeta<\omega$).  By
the assumption on $\lambda$, we may fix a model $\bbM=(\lambda,
\{R^\bbM_{n,\zeta}\}_{n,\zeta <\omega}) $ in the vocabulary $\tau$ with the 
universe $\lambda$ and an ordinal $\alpha^*<\omega_1$ such that: 
\begin{enumerate}
\item[$(\circledast)_{\rm a}$] for every $n$ and a quantifier free formula 
  $\varphi(x_0,\ldots,x_{n-1})\in \cL(\tau)$ there is $\zeta<\omega$ such
  that for all $a_0,\ldots, a_{n-1}\in \lambda$, 
\[\bbM\models\varphi[a_0,\ldots,a_{n-1}]\Leftrightarrow R_{n,\zeta}[a_0,\ldots,
a_{n-1}],\] 
\item[$(\circledast)_{\rm b}$] $\sup\{\rk(v,\bbM):\emptyset\neq v\in
  [\lambda]^{<\omega}\} <\alpha^*$, 
\item[$(\circledast)_{\rm c}$] the rank of every singleton is at least 0. 
\end{enumerate}
For a nonempty finite set $v\subseteq\lambda$ let $\rk(v)=\rk(v,\bbM)$, and 
let  $\zeta(v)<\omega$ and $k(v)<|v|$ be such that $R_{|v|,\zeta(v)},k(v)$ 
witness the rank of $v$. Thus letting $\{a_0,\ldots,a_k,\ldots a_{n-1}\}$ be the
increasing enumeration of $v$ and $k=k(v)$ and $\zeta= \zeta(v)$, we have    
\begin{enumerate}
\item[$(\circledast)_{\rm d}$] if $\rk(v)\geq 0$, then $\bbM\models
  R_{n,\zeta}[a_0,\ldots,a_k,\ldots,a_{n-1}]$ but there is no $a\in
  \lambda\setminus v$  such that  
\[\rk(v\cup\{a\})\geq \rk(v)\ \mbox{ and }\ \bbM\models R_{n,\zeta}
[a_0,\ldots,a_{k-1},a,a_{k+1},\ldots,a_{n-1}],\]  
\item[$(\circledast)_{\rm e}$] if $\rk(v)=-1$, then $\bbM\models R_{n,\zeta}  
  [a_0,\ldots,a_k,\ldots,a_{n-1}]$ but the set 
\[\big\{a\in\lambda:\bbM\models R_{n,\zeta}[a_0,\ldots,a_{k-1},a,a_{k+1},
  \ldots,a_{n-1}]\big\}\] 
is countable.  
\end{enumerate}
Without loss of generality we may also require that (for $\zeta=\zeta(v)$,
$n=|v|$) 
\begin{enumerate}
\item[$(\circledast)_{\rm f}$] for every $b_0,\ldots,b_{n-1}<\lambda$ 
\[\mbox{if }\ \bbM\models R_{n,\zeta}[b_0,\ldots,b_{n-1}] \mbox{ then }\
  b_0<\ldots <b_{n-1}.\]  
\end{enumerate}

\bigskip

Now we will define a forcing notion $\bbP$. {\em A condition\/} $p$ in
$\bbP$ is a tuple 
\[\big(w^p,n^p,\iota^p,M^p,\bar{\eta}^p,\bar{t}^p,\bar{r}^p,\bar{h}^p,
  \bar{g}^p,\cM^p\big)= \big(w,n,\iota,M,\bar{\eta},\bar{t}, \bar{r},
\bar{h},\bar{g},\cM \big)\]     
such that the following demands $(*)_1$--$(*)_{11}$ are satisfied.

\begin{enumerate}
\item[$(*)_1$] $w\in [\lambda]^{<\omega}$, $|w|\geq 5$, $5\leq
  n,M<\omega$, $\iota<\omega$ and if $i^*<\omega$ then $\iota=i^*$.   
\item[$(*)_2$] $\bar{\eta}=\langle \eta_\alpha:\alpha\in w\rangle \subseteq 
  {}^n 2$. 
\item[$(*)_3$] $\bar{t}= \langle t_m:m<M\rangle$, where $\emptyset\neq 
  t_m\subseteq {}^{n\geq} 2$ for $m<M$ is a tree in which all terminal
  branches are of length $n$ and $t_m\cap t_{m'}\cap {}^n 2=\emptyset$ for
  $m<m'<M$. 
\item[$(*)_4$] $\bar{r}=\langle r_m:m<M\rangle$, where $0<r_m\leq n$
  for $m<M$. 
\item[$(*)_5$] $\bar{h}=\langle h_i:i<\iota\rangle$, where $h_i:w^{\langle
    2\rangle}\longrightarrow M$ are such that $h_i(\alpha,\beta)=
  h_i(\beta,\alpha)$. 
\item[$(*)_6$] $\bar{g} =\langle g_i:i<\iota\rangle$, where $g_i:w^{\langle
    2\rangle} \longrightarrow\cO_i$ are such that 
  $\ell\big(g_i(\alpha,\beta)\big)=n$, $g_i(\alpha,\beta)=g_i(\beta,\alpha)$
  and, for each $(\alpha,\beta)\in w^{\langle 2\rangle}$,
  $\big|\bigcup\limits_{i<\iota}  g_i(\alpha,\beta)\big|\geq 6$.
\item[$(*)_7$] For each $m<M$, 
\[t_m\cap {}^n 2=\bigcup\big\{\eta_\alpha+g_i(\alpha,\beta):
  (\alpha,\beta)\in w^{\langle 2\rangle}\mbox{ and }i<\iota\mbox{ and }
  h_i(\alpha,\beta)=m\big\}.\]   
\item[$(*)_8$] The family 
\[\big\{\eta_\alpha:\alpha\in w\big\}\cup\bigcup\big\{g_i(\alpha,\beta):
    (\alpha,\beta)\in w^{\langle 2\rangle}\ \wedge\ i<\iota\big\}\]
is a linearly independent set of vectors in ${}^n2$ (over the field
$\bbZ_2$); in particular there are no repetitions in the representation
above and all elements are non-zero vectors. 
\item[$(*)_9$] $\cM$ consists of all triples $\gd=(\ell^\gd,v^\gd,
  \bm^\gd)=(\ell,v,\bm)$ such that  
\begin{enumerate} 
\item[$(*)_9^{\rm a}$] $0<\ell\leq n$, $v\subseteq w$, $5\leq |v|$,  and
  $\eta_\alpha\rest\ell\neq \eta_\beta\rest\ell$ for distinct $\alpha,\beta\in v$,
\item[$(*)_9^{\rm b}$] $\bm\in \fMtk$, $\ell^\bm=\ell$, $u^\bm=\{\eta_\alpha\rest \ell: \alpha\in v\}$, 
\item[$(*)_9^{\rm c}$] for each $(\alpha,\beta)\in 
  (v)^{\langle 2\rangle}$ and $i<6$ we have 
  $r_{h^\bm_i(\eta_\alpha\rest \ell,\eta_\beta\rest\ell)}\leq \ell^\gd$,
\item[$(*)_9^{\rm d}$] $\big(\forall(\alpha,\beta)\in v^{\langle2\rangle}
  \big) \big(\forall i<6\big) \big(\exists j<\iota\big)
  \big(h^\bm_i(\eta_\alpha \rest \ell,\eta_\beta\rest  \ell)=h_j(\alpha,
  \beta)\big)$.  
\end{enumerate}
\item[$(*)_{10}$] If $\gd_0, \gd_1\in\cM$, $\ell^{\gd_0}=\ell^{\gd_1}=\ell$,
  $\rho\in {}^\ell   2$, and ${\bm^{\gd_1}}={\bm^{\gd_0}}+\rho$, then
  $\rk(v^{\gd_0})=\rk(v^{\gd_1})$, $\zeta(v^{\gd_0})=\zeta(v^{\gd_1})$,
  $k(v^{\gd_0})=k(v^{\gd_0})$ and if $\alpha\in v^{\gd_0}$, $\beta\in
  v^{\gd_1}$ are such that $|\alpha\cap v^{\gd_0}|=k(v^{\gd_0}) =
  k(v^{\gd_1})=|\beta\cap v^{\gd_1}|$,  then $(\eta_\alpha\rest
  \ell)+\rho=\eta_\beta\rest\ell$.    
\item[$(*)_{11}$] Suppose that 
  \begin{itemize}
  \item $\gd_0, \gd_1\in \cM$, $\bm^{\gd_0} \sqsubset \bm^{\gd_1}$ and
    $v^{\gd_0}\subseteq v^{\gd_1}$, and 
\item $\alpha_0\in v^{\gd_0}$,  $|\alpha_0\cap v^{\gd_0}|= k(v^{\gd_0})$,
  $\rk(v^{\gd_0})=-1$. 
  \end{itemize}
Then  $|\{\nu\in u^{\bm^{\gd_1}}:(\eta_{\alpha_0}\rest\ell^{\gd_0})
\trianglelefteq \nu\}|=1$.
\end{enumerate}
To define {\em the order $\leq$ of $\bbP$\/} we declare for $p,q\in \bbP$ that
$p\leq q$\quad if and only if 
\begin{itemize}
\item $w^p\subseteq w^q$, $n^p\leq n^q$, $M^p\leq M^q$, $\iota^p\leq
  \iota^q$ and  
\item $t^p_m=t^q_m\cap {}^{n^p\geq} 2$ and $r^p_m=r^q_m$ for all $m<M^p$,
  and     
\item $\eta^p_\alpha\trianglelefteq \eta^q_\alpha$ for all $\alpha\in w^p$,
  and    
\item $h^q_i\rest (w^p)^{\langle 2\rangle}= h^p_i$ and $g^p_i(\alpha,\beta)
  \preceq_i g^q_i(\alpha,\beta)$ for $i<\iota^p$ and $(\alpha,\beta)\in 
  (w^p)^{\langle 2\rangle}$.
\end{itemize}
\medskip

\begin{claim}
\label{cl8}
\begin{enumerate}
\item $(\bbP,\leq)$ is a partial order of size $\lambda$. 
\item For each $\beta<\lambda$ and $n_0,M_0<\omega$ the set 
\[D_\beta^{n_0,M_0}=\big \{p\in\bbP: n^p>n_0\ \wedge\ M^p>M_0\ \wedge\
\beta\in w^p \big\}\]
is open dense in $\bbP$.
\item If $i^*=\omega$, then for each $\iota<\omega$ the set
  $D_\iota=\{p\in\bbP:\iota^p\geq \iota\}$ is open dense in $\bbP$.
\end{enumerate}
\end{claim}

\begin{proof}[Proof of the Claim]
  (1)\quad First let us argue that $\bbP\neq\emptyset$. Let $\iota=i^*$ if it
  is finite, and $\iota=6$ if $i^*=\omega$. Let $w=\{\alpha_0 , \alpha_1,
  \alpha_2, \alpha_3, \alpha_4\}$ be any 5 element subset of
  $\lambda$. Using \ref{basedef}(2c)+\ref{nicedef}(ii) we may find
  $v(i,b)$ for $i<\iota$ and $b<2$ such that for some
  $\ell<\omega$ for all  $i<\iota$ and $b<2$ we have 
\[v(i,b)\in \cO_i,\quad  v(i,0)\prec_i v(i,1),\quad \mbox{ and
  }\quad \ell\big(v(i,1)\big)=\ell.\]
By \ref{nicedef}(i), we may also require that if $i^*<6$ then for some
$i<\iota$ we have $|v(i,1)|\geq 6$.  Fix an enumeration
\[\big\{(\sigma_a,i_a,j_a,k_a):a<A\big\}=\big\{(\sigma, i, j, k): j<k<5\
  \wedge\ i<\iota\ \wedge\ \sigma\in v(i,1)\big\}.\]
Choose $n>\ell+5$ and a sequence $\langle \rho_a:a<A+5\rangle\subseteq {}^n
2$ so that
\begin{itemize}
\item $\langle \rho_a\rest [\ell,n):a<A+5\rangle$ is linearly independent in
  ${}^{[\ell,n)} 2$ over $\bbZ_2$, and
\item $\sigma_a\vtl\rho_a$ for each $a<A$.   
\end{itemize}
Put
\begin{itemize}
\item $\eta_{\alpha_b}=\rho_{A+b}$ (for $b<5$) and $\bar{\eta}=\langle
  \eta_{\alpha_b}: b<5\rangle$,
\item $g_i(\alpha_j,\alpha_k)=g_i(\alpha_k,\alpha_j)= \big\{\rho_a: a<A\ \wedge\ j=j_a\
  \wedge\ k=k_a\ \wedge\ i_a=i\big\}$ (for $i<\iota$ and $j<k<5$) and
  $\bar{g}=\langle g_i:i<\iota\rangle$.
\end{itemize}
It follows from Definition \ref{nicedef}(iii) that
$g_i(\alpha_j,\alpha_k)\in\cO_i$.

We also let $M=10\cdot \iota$ and we fix a bijection $\varphi:[w]^2 \times
\iota\longrightarrow M$. Then for $j<k<5$ and $i<\iota$ we set
$h_i(\alpha_j,\alpha_k)= h_i(\alpha_k,\alpha_j)=\varphi\big(\{
\alpha_j,\alpha_k\}, i\big)$. This way we defined
$\bar{h}=\langle h_i:i<\iota\rangle$.

We put $r_m=n$ for $m<M$ and we let $t_m\subseteq {}^{n\geq} 2$ be trees in
which all terminal branches are of length $n$ and such that 
\[t_m\cap {}^n2=\bigcup\big\{\eta_\alpha+g_i(\alpha,\beta):
  (\alpha,\beta)\in w^{\langle 2\rangle}\mbox{ and }i<\iota\mbox{ and }
  h_i(\alpha,\beta)=m\big\}.\]
Finally, $\cM$ is defined by clause $(*)_9$.

One easily verifies that $(w,n,\iota,M,\bar{\eta},\bar{t},
\bar{r},\bar{h},\bar{g}, \cM)\in \bbP$.

We see from the arguments above that $|\bbP|\geq\lambda$ and since there are
only countably many elements $p$ of $\bbP$ with $w^p=w$, we get
$|\bbP|=\lambda$.

Clearly, $\leq$ is a partial order on $\bbP$. 
\medskip

\noindent (2)\quad Let $p\in\bbP$, $\beta\in\lambda\setminus w^p$.

We will define a condition $q$ in a manner similar to the construction
in (1) above. Let $\alpha^-=\min(w^p)$ and $\alpha^+=\max(w^p)$.

Set $w^q=w^p\cup\{\beta\}$, $\iota^q=\iota^p$.

For $(\alpha_0,\alpha_1)\in (w^q)^{\langle 2\rangle}$ and $i<\iota^q$
pick $v(i,\alpha_0,\alpha_1)\in \cO_i$ so that: for some $\ell$, for
all $i<\iota^q$ and $(\alpha_0,\alpha_1)\in (w^q)^{\langle 2\rangle}$
we have 
\begin{itemize}
\item $\ell\big(v(i,\alpha_0,\alpha_1)\big)=\ell$,
\item if $\alpha_0,\alpha_1\in w^p$ then $g_i^p(\alpha_0,\alpha_1)
  \prec_i v(i,\alpha_0,\alpha_1) =v(i,\alpha_1,\alpha_0)$,
\item if $\alpha_0\in w^p$ then $g_i^p(\alpha^+,\alpha^-)\prec_i
  v(i,\alpha_0,\beta) =v(i,\beta,\alpha_0)$.
\end{itemize}

Fix an enumeration
\[  \begin{array}{r}
  \big\{(\sigma^a,i^a,\alpha^a_0,\alpha^a_1):a<A\big\}=\big\{(\sigma, 
  i, \alpha_0, \alpha_1): \alpha_0<\alpha_1\mbox{ are from $w^q$ and }\ \ \ \\
  i<\iota^q\ \wedge\ \sigma\in v(i,\alpha_0,\alpha_1)\big\}.
  \end{array}\] 
Choose $n>\ell+|w^p|+1$ and a sequence $\langle \rho_a:a\leq A+|w^p|
\rangle\subseteq {}^n 2$ so that
\begin{itemize}
\item $\langle \rho_a\rest [\ell,n):a\leq A+|w^p|\rangle$ is linearly
  independent in ${}^{[\ell,n)} 2$ over $\bbZ_2$,
\item $\sigma^a\vtl\rho_a$ for each $a<A$, and
\item  if $\alpha\in w^p$ is such that $|w^p\cap \alpha|=k$ then
  $\eta_\alpha^p\vtl \rho_{A+k}$.
\end{itemize}
Put
\begin{itemize}
\item $\eta_\beta^q=\rho_{A+|w^p|}$, and if $\alpha\in w^p$ is such that
  $|w^p\cap \alpha|=k$ then $\eta_\alpha^q=\rho_{A+k}$ and
  $\bar{\eta}^q=\langle \eta^q_\alpha: \alpha\in w^q\rangle$,
\item $g_i^q(\alpha_0,\alpha_1)=g_i^q(\alpha_1,\alpha_0)=
  \big\{\rho_a: a<A\ \wedge\ i=i^a\ \wedge\ \alpha_0=\alpha_0^a\
  \wedge\ \alpha_1=\alpha^a_1\big\}$ (for $i<\iota^q$ and
  $\alpha_0<\alpha_1$ from $w^q$) and $\bar{g}^q=\langle 
  g_i^q:i<\iota^q\rangle$.
\end{itemize}
It follows from Definition \ref{nicedef}(iii) that
$g_i^q(\alpha_0,\alpha_1)\in\cO_i$ and if $(\alpha_0,\alpha_1) \in
(w^p)^{\langle 2\rangle}$ then $g^p_i(\alpha_0,\alpha_1) \prec_i
g^q_i(\alpha_0,\alpha_1)$. 

We also let $M^q=M^p+\iota^q\cdot |w^p|$ and we define mappings 
$h_i^q:(w^q)^{\langle 2\rangle} \longrightarrow M^q$ so that:
\begin{itemize}
\item if $(\alpha_0,\alpha_1)\in (w^p)^{\langle 2\rangle}$ and
  $i<\iota^q$, then $h_i^q(\alpha_0,\alpha_1)=h^p_i(\alpha_0,
  \alpha_1)$,  
\item if $\alpha\in w^p$ and $i<\iota^q$, then
  $h_i^q(\alpha,\beta)=h^q_i(\beta,\alpha)=M^p +|\alpha\cap w^p|\cdot
  \iota+i$.  
\end{itemize}
This way we defined $\bar{h}^q=\langle h_i^q:i<\iota^q\rangle$.

We put $r_m^q=r^p_m$ for $m<M^p$ and $r_m^q=n$ for $M^p\leq m<M^q$. We
let $t_m^q\subseteq {}^{n\geq} 2$ be trees in which all terminal branches
are of length $n$ and such that  
\[t_m^q\cap {}^n2=\bigcup\big\{\eta_\alpha^q+g_i^q(\alpha,\beta):
  (\alpha,\beta)\in (w^q)^{\langle 2\rangle}\mbox{ and
  }i<\iota^q\mbox{ and } h_i^q(\alpha,\beta)=m\big\}.\]
[Note that by our definitions above and by clause $(*)_7$ for $p$ we have 
$t^p_m\cap {}^{n^p}2=t_m^q\cap {}^{n^p}2$ for all $m<M^p$.]
Naturally we also set $n^q=n$ and we define $\cM^q$ by clause $(*)_9$.

We claim that $q=\big(w^q,n^q,\iota^q,M^q,\bar{\eta}^q,\bar{t}^q, \bar{r}^q,
\bar{h}^q,\bar{g}^q, \cM^q\big)\in \bbP$.  Demands $(*)_1$--$(*)_9$ are
pretty straightforward. 
\medskip

\noindent {\bf RE $(*)_{10}$\ :}\quad To justify clause $(*)_{10}$, suppose
that $\gd_0,\gd_1\in  \cM^q$, $\ell^{\gd_0}=\ell^{\gd_1}=\ell$, $\rho\in {}^\ell
2$ and $\bm=\bm^{\gd_0}=\bm^{\gd_1}+\rho$, and consider the following two
cases.   
\medskip

\noindent {\sc Case 1:}\quad $\beta\notin v^{\gd_0}\cup v^{\gd_1}$\\
If $\ell\leq n^p$ then $r_{h_i^\bm(\eta_{\alpha_0}\rest
  \ell,\eta_{\alpha_1}\rest \ell)}\leq n^p$, so  $h_i^\bm(\eta_{\alpha_0}\rest
  \ell,\eta_{\alpha_1}\rest \ell)<M^p$ for all $(\alpha_0,\alpha_1)\in \big(
  v^{\gd_0} \big)^{\langle 2\rangle}$. Hence also $\gd_0,\gd_1\in \cM^p$ and
  clause $(*)_{10}$ for $p$ applies. If $\ell>n^p$ then the sequence
  $\langle \eta^q_\alpha\rest \ell: \alpha\in v^{\gd_0}\cup
  v^{\gd_1}\rangle$ is linearly independent and  
\[\{(\eta^q_\alpha\rest \ell)+\rho:\alpha\in v^{\gd_0}\}=
  \{\eta^q_\alpha\rest \ell:\alpha\in v^{\gd_1}\}.\]
Since $|v^{\gd_0}|\geq 5$ we immediately conclude $\rho={\mathbf 0}$, and 
therefore also $v^{\gd_0}=v^{\gd_1}$ (remember $\ell>n^p$).  
\medskip

\noindent {\sc Case 2:}\quad $\beta\in v^{\gd_0}\cup v^{\gd_1}$\\
Say, $\beta\in v^{\gd_0}$. If $\alpha\in v^{\gd_0}\setminus \{\beta\}$, then
$h^q_j(\alpha,\beta)\geq M^p$ for all $j<\iota$,  and hence
$r^q_{h^\bm_i(\eta_\alpha\rest\ell, \eta_\beta\rest\ell)}=n^q$ (remember
$(*)_9^{\rm d}$). Consequently, $\ell=n^q$. Since the sequence 
$\langle \eta^q_\alpha: \alpha\in v^{\gd_0}\cup v^{\gd_1}\rangle$ is
linearly independent, like before we get $\rho={\mathbf 0}$ and
$v^{\gd_0}=v^{\gd_1}$.  
\medskip

\noindent {\bf RE $(*)_{11}$\ :}\quad Assume towards  contradiction that for 
some $\gd_0,\gd_1\in\cM^q$ we have:
\begin{itemize}
\item $v^\gd_0\subseteq v^\gd_1$ and without loss of generality
  $|v^{\gd_1}|= |v^{\gd_0}|+1$,
\item $\alpha_0\in v^{\gd_0}$, $|\alpha_0\cap v^{\gd_0}|= k(v^{\gd_0})$,
  $\rk(v^{\gd_0})=-1$, and $\bm^{\gd_0}\sqsubset\bm^{\gd_1}$, and 
\item there is $\alpha_1\in  v^{\gd_1}$ such that $\eta^q_{\alpha_0}\rest
  \ell^{\gd_0}= \eta^q_{\alpha_1}\rest  \ell^{\gd_0}$  but
  $\eta^q_{\alpha_0}\rest \ell^{\gd_1}\neq \eta^q_{\alpha_1}\rest \ell^{\gd_1}$.
\end{itemize}
Let $\ell_0=\ell^{\gd_0}$, $\ell_1=\ell^{\gd_1}$.

Suppose $\beta\in v^{\gd_0}$ and take $\beta'\in v^{\gd_0}\setminus
\{\beta\}$. Then $h_j^q(\beta,\beta')\geq M^p$ for all $j<\iota$. Hence, for
some $j<\iota$, 
\[r^q_{h^{\bm^{\gd_0}}_0(\eta_\beta\rest\ell_0,  \eta_{\beta'}\rest\ell_0)}
  =r_{h^p_j(\beta,\beta')}^q=n^q=\ell_0=\ell_1,\]
contradicting the last item in our assumptions.

If we had $v^{\gd_1}=v^{\gd_0}\cup\{\beta\}$, then considering a $\beta'\in 
v^{\gd_0}\setminus \{\alpha_0\}$ we will immediately arrive to  
\[M^p>h^{\bm^{\gd_0}}_0(\eta_\alpha\rest\ell_0,  \eta_{\beta'}\rest\ell_0) =
  h^{\bm^{\gd_1}}_0(\eta_\beta\rest\ell_1,  \eta_{\beta'}\rest\ell_1)\geq
  M^p,\]
a contradiction.

Therefore the only remaining possibility is that $\beta\notin v^{\gd_1}$.

If $\ell_1\leq n^p$, then $\gd_0,\gd_1\in\cM^p$ and clause $(*)_{11}$ for $p$
gives us a contradiction. So assume $\ell_1>n^p$. Since
$\{\eta^q_\gamma\rest n^p:\gamma\in v^{\gd_1}\}$  are all pairwise distinct,
we conclude $\ell_0<n^p$ and $\bm^{\gd_0}\in\cM^p$. We define $\bn\in\fMtk$
by setting: 
\begin{itemize}
\item $\ell^\bn=n^p$, $u^\bn=\{\eta^q_\gamma\rest n^p:\gamma\in v^{\gd_1}\}
  = \{\eta^p_\gamma:\gamma\in v^{\gd_1}\}$, $\iota^\bn=6$,

  and for $(\gamma,\gamma')\in (v^{\gd_1})^{\langle 2\rangle}$ and $i<6$:
\item if $\{\gamma,\gamma'\}\neq \{\alpha_0,\alpha_1\}$, then
\[g^\bn_i(\eta^p_\gamma,\eta^p_{\gamma'})= \{\sigma\rest n^p\!: \sigma\in
  g^{\bm^{\gd_1}}_i(\eta^q_\gamma\rest\ell_1,\eta^q_{\gamma'}\rest\ell_1)
  \}\]  
and $h^\bn_i(\eta^p_\gamma,\eta^p_{\gamma'})=
  h^{\bm^{\gd_1}}_i(\eta^q_\gamma\rest\ell_1,\eta^q_{\gamma'}\rest\ell_1)$,
\item if $\{\gamma,\gamma'\}= \{\alpha_0,\alpha_1\}$, then we fix 
  distinct $\sigma_0,\ldots,\sigma_5\in \bigcup\limits_{j<\iota^q}
  g^p_j(\alpha_0,\alpha_1)$ (remember $(*)_6$ for $p$), and we let
  $g^\bn_i(\eta_{\alpha_0}^p,\eta_{\alpha_1}^p) =g^\bn_i(\eta_{\alpha_1}^p,
  \eta_{\alpha_0}^p) = \{\sigma_i\}$ and
  $h^\bn_i(\eta_{\alpha_0}^p,\eta_{\alpha_1}^p) = h^\bn_i(\eta_{\alpha_1}^p,
  \eta_{\alpha_0}^p)  =m$ where $\eta_{\alpha_0}^p+\sigma_i,
  \eta_{\alpha_1}+\sigma_i\in t^p_m$ (for $i<6$). 
\end{itemize}
Since $\bm^{\gd_0}\sqsubset \bm^{\gd_1}$, in the case when $\{\gamma,
\gamma'\}\neq \{\alpha_0,\alpha_1\}$ we have 
\[g^{\bm^{\gd_0}}_i (\eta^p_\gamma\rest\ell_0, \eta^p_{\gamma'}\rest
  \ell_0)\prec_{\cO^0} g^{\bm^{\gd_1}}_i (\eta^p_\gamma\rest\ell_1,
  \eta^p_{\gamma'}\rest \ell_1)\] 
and hence $g^\bn_i(\eta^p_\gamma, \eta^p_{\gamma'}) \cap
g^\bn_j(\eta^p_\gamma, \eta^p_{\gamma'}) =\emptyset$ whenever $i<j<6$.
Hence \ref{mtkDef}(c) is satisfied. Other cases and other conditions of
\ref{mtkDef} follow immediately by our choices, and hence
\[\bn=(n^p,6,u^\bn,\bar{h}^\bn,\bar{g}^\bn)\in\fMtk.\]
Moreover, $\bm^{\gd_0}\sqsubset \bn$ and
$\gd_*=(n^p,v^{\gd_1},\bn)\in\cM^p$. However, then $\gd_0,\gd_*$ contradict
clause $(*)_{11}$ for $p$. 
\bigskip

\noindent (3)\quad Let $p\in\bbP$. Set $w^q=w^p$ and $\iota^q=\iota^p+1$.
For $(\alpha_0,\alpha_1)\in (w^q)^{\langle 2\rangle}$ and $i<\iota^q$ we use
Proposition \ref{niceprop} to pick $v(i,\alpha_0,\alpha_1)\in \cO_i$ so
that: for some $\ell$, for all $i<\iota^q$ and $(\alpha_0,\alpha_1)\in
(w^q)^{\langle 2\rangle}$ we have 
\begin{itemize}
\item $\ell\big(v(i,\alpha_0,\alpha_1)\big)=\ell$,
\item if $i<\iota^p$ then $g_i^p(\alpha_0,\alpha_1)
  \prec_i v(i,\alpha_0,\alpha_1) =v(i,\alpha_1,\alpha_0)$,
\item for some $v\in\cO_{\iota^p}$, $v\prec_{\iota^p} v(\iota^p,
  \alpha_0,\alpha_1) =v(\iota^p,\alpha_1,\alpha_0)$. 
\end{itemize}

Fix an enumeration
\[  \begin{array}{r}
  \big\{(\sigma^a,i^a,\alpha^a_0,\alpha^a_1):a<A\big\}=\big\{(\sigma, 
  i, \alpha_0, \alpha_1): \alpha_0<\alpha_1\mbox{ are from $w$ and }\ \ \ \\
  i<\iota^q\ \wedge\ \sigma\in v(i,\alpha_0,\alpha_1)\big\}.
  \end{array}\] 
Choose $n=n^q>\ell$ and a sequence $\langle \rho_a:a<A+|w^p|\rangle
\subseteq {}^n 2$ so that 
\begin{itemize}
\item $\langle \rho_a\rest [\ell,n):a<A+|w^q|\rangle$ is linearly
  independent in ${}^{[\ell,n)} 2$ over $\bbZ_2$, and
\item $\sigma^a\vtl\rho_a$ for each $a<A$, and
\item  if $\alpha\in w^q$ is such that $|w^q\cap \alpha|=k$ then
  $\eta_\alpha^p\vtl \rho_{A+k}$.
\end{itemize}
Put
\begin{itemize}
\item if $\alpha\in w^q$ is such that $|w^q\cap \alpha|=k$ then
  $\eta_\alpha^q=\rho_{A+k}$ and 
  $\bar{\eta}^q=\langle \eta^q_\alpha: \alpha\in w^q\rangle$,
\item $g_i^q(\alpha_0,\alpha_1)=g_i^q(\alpha_1,\alpha_0)=
  \big\{\rho_a: a<A\ \wedge\ i=i^a\ \wedge\ \alpha_0=\alpha_0^a\
  \wedge\ \alpha_1=\alpha^a_1\big\}$ (for $i<\iota^q$ and
  $\alpha_0<\alpha_1$ from $w^q$) and $\bar{g}^q=\langle 
  g_i^q:i<\iota^q\rangle$.
\end{itemize}
It follows from Definition \ref{nicedef}(iii) that
$g_i^q(\alpha_0,\alpha_1)\in\cO_i$ and if $(\alpha_0,\alpha_1) \in
(w^p)^{\langle 2\rangle}$ then $g^p_i(\alpha_0,\alpha_1) \prec_i
q^q_i(\alpha_0,\alpha_1)$. 

We also let $M^q=M^p+ |[w^q]^2|$ and we fix a bijection $\psi:
[w^q]^2\longrightarrow [M^p,M^q)$. Then we define mappings 
$h_i^q:(w^q)^{\langle 2\rangle} \longrightarrow M^q$ so that for
$\alpha_0<\alpha_1$ from $w^q$ we have 
\begin{itemize}
\item if $i<\iota^q$, then $h_i^q(\alpha_0,\alpha_1)=
  h_i^q(\alpha_1,\alpha_0)= h^p_i(\alpha_0,\alpha_1)$,  
\item $h_{\iota^p}^q(\alpha_0,\alpha_1)=  h_{\iota^p}^q(\alpha_1,
  \alpha_0)=\psi(\{\alpha_0,\alpha_1\})$.   
\end{itemize}
This way we defined $\bar{h}^q=\langle h_i^q:i<\iota^q\rangle$.

We put $r_m^q=r^p_m$ for $m<M^p$ and $r_m^q=n$ for $M^p\leq m<M^q$. We
let $t_m^q\subseteq {}^{n\geq} 2$ be trees in which all terminal branches
are of length $n$ and such that  
\[t_m^q\cap {}^n2=\bigcup\big\{\eta_\alpha^q+g_i^q(\alpha,\beta):
  (\alpha,\beta)\in (w^q)^{\langle 2\rangle}\mbox{ and
  }i<\iota^q\mbox{ and } h_i^q(\alpha,\beta)=m\big\}.\]
[Note that by our definitions above and by clause $(*)_7$ for $p$ we have
$t^p_m\cap {}^{n^p}2=t_m^q\cap {}^{n^p}2$ for all $m<M^p$.]
We define $\cM^q$ by clause $(*)_9$. Like previously, one easily verifies
that $q=\big(w^q,n^q,\iota^q,M^q,\bar{\eta}^q,\bar{t}^q,\bar{r}^q,
\bar{h}^q,\bar{g}^q, \cM^q\big)\in \bbP$. 
[The crucial point is that if $\gd\in\cM^q$, $\eta,\nu\in u^{\bm^\gd}$ and
$h^{\bm^\gd}_i(\eta,\nu) \geq M^p$, then $\ell^\gd=n^q$.]
\end{proof}

\begin{claim}
  \label{cl9}
The forcing notion $\bbP$ has the Knaster property.
\end{claim}

\begin{proof}[Proof of the Claim]
Suppose that $\langle p_\xi:\xi<\omega_1\rangle$ is a sequence of pairwise 
distinct conditions from $\bbP$ and let
\[p_\xi=\big(w_\xi,n_\xi,\iota_\xi, M_\xi,\bar{\eta}_\xi, \bar{t}_\xi,
  \bar{r}_\xi,  \bar{h}_\xi,\bar{g}_\xi, \cM_\xi\big)\]
where $\bar{\eta}_\xi=\langle\eta^\xi_\alpha: \alpha\in w_\xi\rangle$,
$\bar{t}_\xi=\langle t^\xi_m:m<M_\xi\rangle$, $\bar{r}_\xi=\langle
r^\xi_m:m<M_\xi\rangle$,  and $\bar{h}_\xi= \langle
h^\xi_i:i<\iota_\xi\rangle$, $\bar{g}_\xi=\langle g^\xi_i:i_\xi<\iota\rangle$.  By a
standard $\Delta$--system cleaning procedure we may find an uncountable set
$A\subseteq \omega_1$ such that the following demands $(*)_{12}$--$(*)_{15}$
are satisfied. 
\begin{enumerate}
\item[$(*)_{12}$] $\{w_\xi:\xi\in A\}$ forms a $\Delta$--system with the
  kernel $w^*$. 
\item[$(*)_{13}$] If  $\xi,\varsigma\in A$, then $|w_\xi|=|w_{\varsigma}|$ ,
  $n_\xi=n_{\varsigma}$, $\iota_\xi=\iota_\varsigma$, $M_\xi=M_{\varsigma}$,  and
  $t^\xi_m=t^{\varsigma}_m$ and $r^\xi_m=r^\varsigma_m$ (for  $m<M_\xi$).      
\item[$(*)_{14}$] If $\xi<\varsigma$ are from $A$ and
  $\pi:w_\xi\longrightarrow w_{\varsigma}$ is the order isomorphism, then 
  \begin{enumerate}
  \item[(a)] $\pi(\alpha)=\alpha$ for $\alpha\in w^*=w_\xi\cap w_{\varsigma}$, 
  \item[(b)] if $\emptyset\neq v\subseteq w_\xi$, then $\rk(v)=\rk(\pi[v])$,
    $\zeta(v)= \zeta(\pi[v])$ and $k(v)=k(\pi[v])$, 
  \item[(c)] $\eta_\alpha^\xi=\eta_{\pi(\alpha)}^{\varsigma}$ (for
    $\alpha\in  w_\xi$), 
  \item[(d)] $g_i^\xi(\alpha,\beta)=g_i^\zeta(\pi(\alpha),\pi(\beta))$ and
    $h_i^\xi(\alpha,\beta) =h_i^\zeta(\pi(\alpha),\pi(\beta))$ for
    $(\alpha,\beta)\in (w_\xi)^{\langle 2\rangle}$ and $i<\iota_\xi$, and 
  \end{enumerate}
\item[$(*)_{15}$] $\cM_\xi=\cM_\varsigma$ (this actually follows from the
previous demands). 
\end{enumerate}
Note that then also 
\begin{enumerate}
\item[$(*)_{16}$] if $\xi\in A$, $v\subseteq w^*$ and $\delta\in
  w_\xi\setminus w^*$ are such that $\rk\big(v\cup\{\delta\}\big)=-1$, then 
  $k\big(v\cup\{\delta\}\big)\neq |\delta\cap v|$.
\end{enumerate}
[Why? Suppose $\rk\big(v\cup\{\delta\}\big)=-1$ and
$k=k\big(v\cup\{\delta\}\big)= |\delta\cap v|$, $j=j\big( v\cup
\{\delta\}\big)$.  For $\varsigma\in A$ let  $\pi_\varsigma: w_\xi
\longrightarrow w_\varsigma$ be the order isomorphism and let
$\delta_\varsigma=\pi_\varsigma(\delta)$. By $(*)_{14}$ we know that
$k=k\big(v\cup\{\delta_\varsigma\}\big) =|\delta_\varsigma\cap v|$ and
$j=j\big(v\cup\{\delta_\varsigma\}\big)$. Therefore, letting
$v\cup\{\delta\}=\{a_0,\ldots,a_{n-1}\}$ be the increasing enumeration, for
every $\varsigma\in A$ we have $\bbM\models R_{n,j}[a_0,\ldots,a_{k-1},
\delta_\varsigma, a_{k+1},\ldots, a_{n-1}]$. Hence the set 
\[\{b<\lambda: \bbM\models R_{n,j}[a_0,\ldots,a_{k_1}, b , a_{k+1},\ldots,
  a_{n-1}] \}\]
is uncountable, contradicting $(\circledast)_{\rm e}$ from the beginning of
the proof of the theorem.]
\medskip

We will show that for distinct $\xi,\varsigma$ from $A$ the conditions
$p_\xi,p_\varsigma$ are compatible. So let $\xi,\varsigma\in A$,
$\xi<\varsigma$ and let $\pi:w_\xi\longrightarrow w_\varsigma$ be the order
isomorphism. We will define $q=\big(w,n,\iota, M,\bar{\eta}, \bar{t},
\bar{r}, \bar{h},\bar{g},  \cM\big)$ where $\bar{\eta}=\langle\eta_\alpha:
\alpha\in w\rangle$, $\bar{t}=\langle t_m:m<M\rangle$, $\bar{r}=\langle
r_m:m<M\rangle$, and $\bar{h}= \langle h_i:i<\iota \rangle$,
$\bar{g}=\langle g_i:i<\iota\rangle$.   
\medskip

We set
\begin{enumerate}
\item[$(*)_{17}$] $\iota=\iota_\xi$ and $w=w_\xi\cup w_\varsigma$.
\end{enumerate}
\medskip

Similarly to the arguments in previous claims, we first pick
\[\langle v(i,\alpha_0,\alpha_1): (\alpha_0,\alpha_1)\in w^{\langle
    2\rangle}\ \wedge\ i<\iota\rangle\]
and an $\ell$ such that for all $i<\iota$ and $(\alpha_0,\alpha_1)\in
w^{\langle 2\rangle}$ we have  
\begin{itemize}
\item $v(i,\alpha_0,\alpha_1) =v(i,\alpha_1,\alpha_0)\in \cO_i$,
  $\ell\big(v(i,\alpha_0,\alpha_1)\big)=\ell$, 
\item if $\alpha_0,\alpha_1\in w_\xi$ then $g_i^\xi(\alpha_0,\alpha_1)
  \prec_i v(i,\alpha_0,\alpha_1) $, and 
\item  if $\alpha_0,\alpha_1\in w_\varsigma$ then $g_i^\varsigma
  (\alpha_0,\alpha_1) \prec_i v(i,\alpha_0,\alpha_1)$.
\end{itemize}
(No other demands on $v(i,\alpha_0,\alpha_1)$ but symmetry if $\alpha_0\in 
w_\xi\setminus w_\varsigma$ and $\alpha_1\in w_\varsigma\setminus w_\xi$.)
Then we fix an enumeration
\[  \begin{array}{r}
  \big\{(\sigma^a,i^a,\alpha^a_0,\alpha^a_1):a<A\big\}=\big\{(\sigma, 
  i, \alpha_0, \alpha_1): \alpha_0<\alpha_1\mbox{ are from $w$ and }\ \ \ \\
  i<\iota\ \wedge\ \sigma\in v(i,\alpha_0,\alpha_1)\big\}
  \end{array}\] 
and we choose $n>\ell$ and $\langle \rho_a:a< A+|w|\rangle\subseteq {}^n 2$ 
so that 
\begin{itemize}
\item $\langle \rho_a\rest [\ell,n):a<A+|w|\rangle$ is linearly
  independent in ${}^{[\ell,n)} 2$ over $\bbZ_2$, and
\item $\sigma^a\vtl\rho_a$ for each $a<A$, and
\item  if $\alpha\in w_\xi$ is such that $|w_\xi\cap \alpha|=k$ then
  $\eta_\alpha^\xi\vtl \rho_{A+k}$,
\item if $\alpha\in w_\varsigma\setminus w_\xi$ is such that
  $|(w_\varsigma\setminus w_\xi)\cap \alpha|=k$ then $\eta_\alpha^\varsigma
  \vtl \rho_{A+|w_\xi|+k}$. 
\end{itemize}
Put

\begin{enumerate}
\item[$(*)_{18}$] $n$ is the one chosen right above,
\item[$(*)_{19}$] $\bar{\eta}=\langle\eta_\alpha:\alpha\in w\rangle$, where  
  \begin{itemize}
\item  if $\alpha\in w_\xi$ is such that $|w_\xi\cap \alpha|=k$ then 
  $\eta_\alpha=\rho_{A+k}$, 
\item if $\alpha\in w_\varsigma\setminus w_\xi$ is such that $|(w_\varsigma
  \setminus  w_\xi)\cap \alpha|=k$ then $\eta_\alpha=\rho_{A+|w_\xi|+k}$, 
\end{itemize}
\item[$(*)_{20}$] $\bar{g}=\langle  g_i:i<\iota\rangle$, where for $i<\iota$ and 
  $\alpha_0<\alpha_1$ from $w$ we put 
\[g_i(\alpha_0,\alpha_1)=g_i(\alpha_1,\alpha_0)=
  \big\{\rho_a: a<A\ \wedge\ i=i^a\ \wedge\ \alpha_0=\alpha_0^a\
  \wedge\ \alpha_1=\alpha^a_1\big\}.\]
\end{enumerate}
As before, by \ref{nicedef}(iii), we know that 
$g_i(\alpha_0,\alpha_1)\in\cO_i$ and if $(\alpha_0,\alpha_1) \in
(w_\xi)^{\langle 2\rangle}$ then $g^\xi_i(\alpha_0,\alpha_1) \prec_i
g_i(\alpha_0,\alpha_1)$ and similarly for $\varsigma$ in place of $\xi$.   

Let
\begin{enumerate}
\item[$(*)_{21}$]  $M=M_\xi+ |w_\xi\setminus w_\varsigma|^2$
\end{enumerate}
and let $\psi: (w_\xi\setminus w_\varsigma)\times (w_\varsigma \setminus 
w_\xi)\longrightarrow [M_\xi,M)$ be a bijection. Then we define
\begin{enumerate}
\item[$(*)_{22}$] $\bar{h}=\langle h_i:i<\iota\rangle$, where mappings  
$h_i:w^{\langle 2\rangle} \longrightarrow M$ are such that for distinct 
$\alpha_0,\alpha_1\in w$ and $i<\iota$ we have 
\begin{itemize}
\item $h_i(\alpha_0,\alpha_1)=   h_i(\alpha_1,\alpha_0)$,  
\item if $\alpha_0,\alpha_1\in w_\xi$, then $h_i(\alpha_1,\alpha_0)=
  h_i^\xi(\alpha_1,\alpha_0)$, 
\item if $\alpha_0,\alpha_1\in w_\varsigma$, then $h_i(\alpha_1,\alpha_0)=
  h_i^\varsigma (\alpha_1,\alpha_0)$, 
\item if $\alpha_0\in w_\xi\setminus w_\varsigma$ and $\alpha_1\in 
  w_\varsigma\setminus w_\xi$, then $h_i(\alpha_1,\alpha_0)= 
  \psi(\alpha_0,\alpha_1)$. 
\end{itemize}
\item[$(*)_{23}$] $\bar{t}=\langle t_m:m<M\rangle$, where $t_m\subseteq
  {}^{n\geq} 2$ are trees in which all terminal branches are of length $n$ and such that  
\[t_m\cap {}^n2=\bigcup\big\{\eta_\alpha+g_i(\alpha,\beta): 
  (\alpha,\beta)\in w^{\langle 2\rangle}\mbox{ and }i<\iota\mbox{ and }
  h_i(\alpha,\beta)=m\big\},\] 
\item[$(*)_{24}$] $\bar{r}=\langle r_m:m<M\rangle$, where $r_m=r^\xi_m$ for
  $m<M_\xi$, $r_m= n$ if $M_\xi\leq m<M$.
 \item[$(*)_{25}$] $\cM$ is defined by $(*)_9$ (for the objects introduced in
   $(*)_{17}$--$(*)_{24}$).  
\end{enumerate}

In clauses $(*)_{17}$--$(*)_{25}$ we defined all the ingredients of 
\[q=\big(w,n,M,\bar{\eta},\bar{t},\bar{r},\bar{h},\bar{g},\cM\big).\]
We still need to argue that $q\in \bbP$ (after this it will be obvious that
it is a condition stronger than both $p_\xi$ and $p_\varsigma$).

It is  pretty straightforward that $q$ satisfies demands $(*)_1$--$(*)_9$. 
\bigskip

\noindent {\bf RE $(*)_{10}$\ :}\quad To justify clause $(*)_{10}$, suppose
that $\gd_0,\gd_1\in  \cM$, $\ell^{\gd_0}=\ell^{\gd_1}=\ell$ and $\rho\in
{}^\ell 2$ and $\bm^{\gd_1}=\bm^{\gd_0}+\rho$, and consider the following
three cases. 
\medskip

\noindent {\sc Case 1:}\quad $v^{\gd_0}\subseteq w_\xi$\\
Then for each $(\delta,\vare)\in (v^{\gd_0})^{\langle 2\rangle}$ and
$i<\iota$ we have $h_i(\delta,\vare)<M_\xi$, and consequently
$\rng\big(h_j^{\bm^{\gd_0}}\big)\subseteq M_\xi$ (for $j<6$). Hence also 
$\rng\big(h_j^{\bm^{\gd_1}}\big)\subseteq M_\xi$ (for $j<6$). But looking at
$(*)_{22}$ (and remembering $(*)_9^{\rm d}$) we now conclude
$h_i(\delta,\vare)<M_\xi$ for $(\delta,\vare)\in (v^{\gd_1})^{\langle
  2\rangle}$ and $i<\iota$. Consequently, either $v^{\gd_1}\subseteq w_\xi$
or $v^{\gd_1}\subseteq w_\varsigma$. 

If $v^{\gd_1}\subseteq w_\xi$ and $\ell\leq n_\xi$, then $\gd_0,\gd_1\in
\cM_\xi$ and clause $(*)_{10}$ for $p_\xi$ can be used to get the desired
conclusion. 

If $v^{\gd_1}\subseteq w_\xi$ and $\ell>n_\xi$, then $\{\eta_\alpha\rest
\ell: \alpha\in v^{\gd_0}\cup v^{\gd_1}\}$ is linearly independent and hence
$\rho={\mathbf 0}$ and $v^{\gd_0}=v^{\gd_1}$. 

If $v^{\gd_1}\subseteq w_\varsigma$ and $\ell\leq n_\xi$, then consider
$v=\pi^{-1}[v^{\gd_1}]\subseteq w_\xi$ and
$\gd=(\ell,v,\bm^{\gd_1})$. Clearly, $\gd\in\cM_\xi$ and we may use
$(*)_{10}$ for $p_\xi$ to conclude that $\rk(v)=\rk(v^{\gd_0})$,
$\zeta(v)=\zeta(v^{\gd_0})$, $k(v)= k(v^{\gd_0})$, and if $\alpha\in
v^{\gd_0}$, $\beta\in v$ are such that $|\alpha\cap v^{\gd_0}|= k(v^{\gd_0})
=k(v)=|\beta\cap v|$, then $(\eta_\alpha\rest
\ell)+\rho=\eta_\beta\rest\ell$. Now we use the properties $(*)_{14}$(b,c)
of $\pi$ to get a similar assertions with $v^{\gd_1}$ in place of $v$. 

If $v^{\gd_1}\subseteq w_\varsigma$ and $\ell> n_\xi$, then we consider
$v=\pi^{-1}[v^{\gd_1}]\subseteq w_\xi$ and use the linear independence of
$\{\eta_\alpha\rest \ell: \alpha\in v^{\gd_0}\cup v\}$ to conclude that
$\rho={\mathbf 0}$ and $v^{\gd_0}=v=\pi^{-1}\big[v^{\gd_1}\big]$.  Finally we
use the properties $(*)_{14}$(b,c) of $\pi$ to get the desired assertions.
\medskip

\noindent {\sc Case 2:}\quad $v^{\gd_0}\subseteq w_\varsigma$\\
Same as the previous case, just interchanging $\xi$ and $\varsigma$. 
\medskip

\noindent {\sc Case 3:}\quad $v^{\gd_0}\setminus w_\xi\neq \emptyset \neq
v^{\gd_0}\setminus w_\varsigma$\\
Then for some $(\delta,\vare)\in (v^{\gd_0})^{\langle 2\rangle}$ we have
$h_i(\delta,\vare)\geq M_\xi$ for all $i<\iota$, so necessarily
$\ell=n$. Now, the linear independence of $\bar{\eta}$ implies
$\rho={\mathbf 0}$ and $v^{\gd_0}=v^{\gd_1}$ and the desired conclusion
follows.  
\medskip

\noindent {\bf RE $(*)_{11}$\ :}\quad Let us prove clause $(*)_{11}$
now. Suppose that $\gd_0,\gd_1 \in\cM$, $\delta\in v^{\gd_0}$,
$|\delta\cap v^{\gd_0}|= k(v^{\gd_0})$, $\rk(v^{\gd_0})=-1$, and
$v^{\gd_0}\subseteq v^{\gd_1}$ and $\bm^{\gd_0} \sqsubset
\bm^{\gd_1}$. Assume towards contradiction that there is an 
$\vare\in  v^{\gd_1}$ such that 
\begin{enumerate}
\item[$(*)_{26}$] $\eta_\vare\rest \ell^{\gd_1}\neq \eta_\delta \rest
  \ell^{\gd_1}$ but $\eta_\vare\rest
  \ell^{\gd_0}=\eta_\delta\rest\ell^{\gd_0}$.
\end{enumerate}
Without loss of generality $v^{\gd_1}= v^{\gd_0}\cup\{\vare\}$. Since we
must have $\ell^{\gd_0}<n$, for no $\alpha,\beta\in v^{\gd_0}$ we can have
$(\forall i<\iota)(h_i(\alpha,\beta)\geq M_\xi)$. Therefore either 
$v^{\gd_0} \subseteq w_\xi$ or $v^{\gd_0}\subseteq w_\varsigma$. By the
symmetry, we may assume $v^{\gd_0}\subseteq w_\xi$. Note that 
\begin{enumerate}
\item[$(*)_{27}$] if $(\alpha,\beta)\in (v^{\gd_1})^{\langle 
  2\rangle}\setminus \{(\vare,\delta),(\delta, \vare)\}$ then 
$h_i(\alpha,\beta)<M_\xi$ for all $i<\iota$.  
\end{enumerate}
Now, if $v^{\gd_1}\subseteq w_\xi$ and $\ell^{\gd_1}\leq
n_\xi$, then $\gd_0,\gd_1\in\cM_\xi$ and they contradict clause $(*)_{11}$
for $p_\xi$. Let us consider the possibility that $v^{\gd_1}\subseteq
w_\xi$ but $\ell^{\gd_1}>n_\xi$. Define $\bn\in\fMtk$ by:
\begin{itemize}
\item $\ell^\bn=n_\xi$, $u^\bn=\{\eta_\gamma\rest n_\xi:\gamma\in
  v^{\gd_1}\}$ (note $\eta_\vare\rest n_\xi\neq \eta_\delta\rest n_\xi$),
  $\iota^\bn=6$, and for $(\gamma,\gamma')\in  (v^{\gd_1})^{\langle
    2\rangle}$ and $i<6$:  
\item if $\{\gamma,\gamma'\}\neq \{\vare,\delta\}$, then
\[g^\bn_i(\eta_\gamma\rest n_\xi,\eta_{\gamma'}\rest n_\xi)=
  \{\sigma\rest n_\xi: \sigma\in  g^{\bm^{\gd_1}}_i(\eta_\gamma
  \rest \ell^{\gd_1},\eta_{\gamma'} \rest \ell^{\gd_1}) \}\]
and $h^\bn_i(\eta_\gamma\rest n_\xi,\eta_{\gamma'}\rest n_\xi)= 
  h^{\bm^{\gd_1}}_i(\eta_\gamma\rest \ell^{\gd_1},\eta_{\gamma'}\rest
  \ell^{\gd_1})$, and 
\item for $\{\gamma,\gamma'\}= \{\delta,\vare\}$  we fix any
  distinct $\sigma_0,\ldots,\sigma_5\in \bigcup\limits_{j<\iota}
  g_j^\xi(\delta,\vare)$ and we let $g^\bn_i(\eta_\delta\rest n_\xi,
  \eta_\vare\rest n_\xi) =g^\bn_i(\eta_\vare\rest n_\xi,
  \eta_\delta\rest n_\xi) = \{\sigma_i\}$ and $h^\bn_i(\eta_\delta\rest
  n_\xi,\eta_\vare\rest n_\xi) = h^\bn_i(\eta_\vare\rest n_\xi, 
  \eta_\delta\rest n_\xi)  =m$ where $(\eta_\delta\rest n_\xi)+ \sigma_i,  
  (\eta_\vare\rest n_\xi)+\sigma_i\in t^\xi_m$ (for $i<6$).  
\end{itemize}
Since $\bm^{\gd_0}\sqsubset \bm^{\gd_1}$, in the case when $\{\gamma,
\gamma'\}\neq \{\delta,\vare\}$ we have
\[g^{\bm^{\gd_0}}_i (\eta_\gamma\rest\ell^{\gd_0}, \eta_{\gamma'}\rest
  \ell^{\gd_0})\prec_{\cO^0} g^{\bm^{\gd_1}}_i
  (\eta_\gamma\rest\ell^{\gd_1}, \eta_{\gamma'}\rest \ell^{\gd_1}),\] 
and hence $g^\bn_i(\eta^p_\gamma, \eta^p_{\gamma'}) \cap
g^\bn_j(\eta^p_\gamma, \eta^p_{\gamma'}) =\emptyset$ whenever $i<j<6$.
Hence \ref{mtkDef}(c) is satisfied. Other cases and other conditions of
\ref{mtkDef} follow immediately by our choices, and hence
\[\bn=(\ell^\bn,6,u^\bn,\bar{h}^\bn,\bar{g}^\bn)\in\fMtk.\]
Moreover, $\bm^{\gd_0}\sqsubset \bn$ and
$\gd_*=(n_\xi,v^{\gd_1},\bn)\in\cM_\xi$. However, then $\gd_0,\gd_*$ contradict
clause $(*)_{11}$ for $p_\xi$.

Consequently, $v^{\gd_1}\setminus w_\xi\neq \emptyset$, so necessarily
$\vare\notin w^*$.

Suppose $|v^{\gd_0}\setminus w^*|\geq 2$, say $\alpha_0,\alpha_1\in
v^{\gd_0}\setminus w^*$. Then $h_i(\vare,\alpha_0), h_i(\vare,\alpha_1)\geq
M_\xi$ for all $i<\iota$. But $\bm^{\gd_0}\sqsubset \bm^{\gd_1}$ implies
that for $\alpha\in v^{\gd_0}\setminus \{\delta\}$ we have
\[h_0^{\bm^{\gd_1}}(\eta_\vare\rest\ell^{\gd_1},
  \eta_\alpha\rest\ell^{\gd_1}) = h_0^{\bm^{\gd_0}}(\eta_\delta
  \rest\ell^{\gd_0},  \eta_\alpha\rest\ell^{\gd_0})<M_\xi,\]
so we arrive to a contradiction.
\medskip

If we had $v^{\gd_0}\subseteq w^*$, then $v^{\gd_1}\subseteq w_\varsigma$
and we may repeat the earlier arguments with $\varsigma$ in place 
of $\xi$ to get a contradiction. Thus the only possibility left is that 
$|v^{\gd_0}\setminus w^*|=1$. Let $\{\alpha\}=v^{\gd_0}\setminus w^*$. If
$\alpha\neq\delta$, then $h^{\bm^{\gd_1}}_0(\eta_\alpha \rest \ell^{\gd_1},
\eta_\vare\rest \ell^{\gd_1})=h^{\bm^{\gd_0}}_0(\eta_\alpha \rest \ell^{\gd_0},
\eta_\vare\rest \ell^{\gd_0})<M_\xi$ gives a contradiction like
before. Therefore,  $v^{\gd_0}=(v^{\gd_0}\cap w^*)\cup\{\delta\}$. But now
our assumptions on $v^{\gd_0},\delta$ contradict $(*)_{16}$.
\end{proof}

\begin{claim}
\label{cl7}
Assume $p=\big(w,n,\iota,M,\bar{\eta},\bar{t},\bar{h},\bar{g},\cM
\big)\in\bbP$. If $\bm\in \fMtk$ is such that $\ell^\bm=n$ and $|u^\bm|\geq
5$, then for some $\rho\in {}^n 2$ and $v\subseteq w$ we have
$\big(n,v,(\bm+\rho)\big) \in\cM$.
\end{claim}

\begin{proof}[Proof of the Claim]
Let $\bm\in \fMtk$ be such that $\ell^\bm=n$. Suppose $(\eta,\nu)\in
\big(u^\bm\big)^{\langle 2\rangle}$.  

Let $g^\bm_j(\eta,\nu)= \{\sigma_j\}$ for $j<6$. Then $\sigma_j$s are
pairwise distinct, and if $\eta+\sigma_i=\nu+\sigma_j$ then
\[\eta+\sigma_k,\nu+\sigma_k\notin \{\eta+\sigma_i,\nu+\sigma_i\}=
  \{\eta+\sigma_j, \nu+\sigma_j\}\] 
whenever $k\notin \{i,j\}$. Hence we may pick $j_0<j_1<j_2<6$ such that 
\[\eta+\sigma_{j_0}, \nu+\sigma_{j_0}, \eta+\sigma_{j_1}, \nu+\sigma_{j_1},
  \eta+\sigma_{j_2}, \nu+\sigma_{j_2}\]
are all pairwise distinct.  Just to simplify notation let us assume that
$j_0=0$, $j_1=1$ ad $j_2=2$.

For each $j<3$ we have $\eta+\sigma_j, \nu+\sigma_j\in \bigcup_{m<M}
t_m$. By clause $(*)_7$ there are $(\alpha_j,\beta_j), 
(\alpha_j',\beta_j') \in w^{\langle 2\rangle}$ and $\rho_j\in
\bigcup\limits_{i<\iota} g_i(\alpha_j,\beta_j)$ and $\rho_j'\in
\bigcup\limits_{i<\iota} g_i(\alpha_j',\beta_j')$ such that
$\eta+\sigma_j=\eta_{\alpha_j}+\rho_j$ and $\nu+\sigma_j=
  \eta_{\alpha_j'} +\rho_j'$ for $j<3$.Then
  $\eta+\nu=\eta_{\alpha_j}+\eta_{\alpha_j'}+\rho_j+\rho_j'$ for all 
$j<3$. We will consider 3 cases, the first two of them will be shown to be
impossible. 
\medskip

\noindent {\sc Case 1:}\quad $\eta_{\alpha_j}=\eta_{\alpha_j'}$ for some 
$j<3$. \\
Then, by the linear independence demanded in $(*)_7$,
$\eta_{\alpha_j}=\eta_{\alpha_j'}$ for all $j<3$ and $\{\rho_0,\rho_0'\}
=\{\rho_1,\rho_1'\}=\{\rho_2,\rho_2'\}$. But $g_i(\alpha,\beta)$'s are
disjoint, so each $\rho\in \bigcup\big\{ 
g_i(\alpha,\beta):(\alpha,\beta)\in w^{\langle 2\rangle}\ \wedge\
i<\iota\big\}$  uniquely determines $\alpha,\beta$ such that
$\eta_\alpha+\rho,\eta_\beta+\rho\in \bigcup\limits_{m<M} t_m$. Therefore,
$|\{\alpha_0,\alpha_1,\alpha_2\}|\leq 2$ in the current case. Since
$\eta+\sigma_j,\nu+\sigma_j$ are all pairwise distinct (for $j<3$), this
gives an immediate contradiction.  
\medskip

\noindent {\sc Case 2:}\quad $\eta_{\alpha_j}\neq\eta_{\alpha_j'}$ and
$\rho_j\neq \rho_j'$ for some (equivalently: all) $j<3$. \\
Then $\{\eta_{\alpha_0},\eta_{\alpha_0'}\}=
\{\eta_{\alpha_1},\eta_{\alpha_1'}\} = \{\eta_{\alpha_2},\eta_{\alpha_2'}\}$
and $\{\rho_0,\rho_0'\}= \{\rho_1,\rho_1'\} =\{\rho_2,\rho_2'\}$. However,
this again contradicts $\eta+\sigma_j,\nu+\sigma_j$ being pairwise
distinct. 
\medskip

Thus the only possible case is the following:

\noindent {\sc Case 3:}\quad $\eta_{\alpha_j}\neq\eta_{\alpha_j'}$ and
$\rho_j=\rho_j'$ for all $j<3$.\\
Then $\eta+\nu=\eta_{\alpha_0}+\eta_{\alpha_0'}$.
\smallskip

Consequently we have shown that 
\[u^\bm+ u^\bm\subseteq \{\eta_\alpha+\eta_\beta:\alpha,\beta\in w\}.\]
By Lemma \ref{litlem} for some $\rho$ we have $u^\bm+\rho\subseteq
\{\eta_\alpha:\alpha\in w\}$. Let $v=\{\alpha\in w: \eta_\alpha\in
u^\bm+\rho\}$. Let us argue that $\big(n,v,(\bm+\rho)\big)\in \cM$: demands
$(*)^{\rm a}_9$--$(*)^{\rm c}_9$ are immediate consequences of our choices
above. Let us verify $(*)^{\rm d}_9$.

Suppose that $(\alpha,\beta)\in v^{\langle 2\rangle}$ and $i<6$. Let
$\eta=\eta_\alpha+\rho, \nu=\eta_\beta+\rho$ (so they are in $u^\bm$) and
let $\{\sigma_i\}= g^\bm_i(\eta,\nu)$. Then $\eta+\sigma_i, \nu+\sigma_i\in
\bigcup\limits_{m<M} t_m$, so we may choose $(\alpha',\beta'),
(\alpha'',\beta'')\in w^{\langle 2\rangle}$ and $j',j''<\iota$ and
$\rho'\in g_{j'}(\alpha',\beta')$ and $\rho''\in g_{j''}(\alpha'',\beta'')$
such that $\eta+\sigma_i=\eta_{\alpha'}+\rho'$ and  
$\nu+\sigma_i=\eta_{\alpha''}+\rho''$. Then
\[\eta_\alpha+\eta_\beta= \eta+\nu=\eta_{\alpha'}+\eta_{\alpha''}+
  \rho'+\rho''.\]
By the linear independence stated in $(*)_8$ we get $\rho'=\rho''$ and
$\{\eta_{\alpha'},\eta_{\alpha''}\}=\{\eta_\alpha,\eta_\beta\}$. Consequently
also $\{\alpha,\beta\}=\{\alpha',\alpha''\}$ and $\{\alpha',\beta'\}=
\{\alpha'',\beta''\}$ and $j'=j''$.  Since $\alpha\neq\beta$ we get
$\alpha'\neq\alpha''$ and thus $\alpha'=\beta''$,
$\alpha''=\beta'$. Consequently, $\{\alpha'',\beta''\}=\{\alpha',\beta'\}=
\{\alpha',\alpha''\}=\{\alpha,\beta\}$. Hence
$\eta+\sigma_i=\eta_{\alpha'}+\rho'\in 
t_{h_{j'}(\alpha,\beta)}= t_{h_{j'}(\beta,\alpha)}$ and
$\nu+\sigma_i=\eta_{\alpha''}+\rho'\in t_{h_{j'}(\alpha,\beta)}=
t_{h_{j'}(\beta,\alpha)}$. Therefore, 
$h^{\bm+\rho}_i(\eta_\alpha,\eta_\beta)=h^\bm_i(\eta,\nu)=h_{j'}(\alpha,\beta)
  =h_{j'}(\beta,\alpha)$.
\end{proof}
\bigskip

Define $\bbP$--names $\name{T}_m$ and $\name{\eta}_\alpha$ (for
$m<\omega$ and $\alpha<\lambda$) by

$\forces_\bbP$`` $\name{T}_m= \bigcup\{t^p_m: p\in
  \name{G}_\bbP\ \wedge\ m<M^p\}$ '', and 

$\forces_\bbP$`` $\name{\eta}_\alpha= \bigcup\{\eta^p_\alpha:
  p\in \name{G}_\bbP\ \wedge\ \alpha\in w^p\}$ ''.

\begin{claim}
  \label{cl10}
  \begin{enumerate}
\item For each $m<\omega$ and $\alpha<\lambda$,

$\forces_\bbP$`` $\name{\eta}_\alpha\in\can$ and
$\name{T}_m\subseteq {}^{\omega>}2$ is a tree without terminal nodes ''.
\item For all $\alpha<\beta<\lambda$ we have
 \[\forces_\bbP\mbox{`` }\big(\eta_\alpha+\bigcup\limits_{m<\omega}
  \lim(\name{T}_m)\big)\cap  \big(\eta_\beta+\bigcup\limits_{m<\omega} 
  \lim(\name{T}_m)\big)\mbox{ is $\bar{\cO}$--large ''.}\]
\item $\forces_{\bbP}$`` $\bigcup\limits_{m<\omega}
  \lim(\name{T}_m)$ is  a $\bar{\cO}^6$--{\bf npots} set ''.  
  \end{enumerate}
\end{claim}

\begin{proof}[Proof of the Claim]
(1, 2) \quad By Claim \ref{cl8} (and the definition of the order in $\bbP$). 
\medskip

\noindent (3)\quad Let $G\subseteq \bbP$ be a generic filter over $\bV$ and 
let us work in $\bV[G]$.  Let $\bar{T}=\langle
(\name{T}_m)^G:m<\omega\rangle$.  

Suppose towards contradiction that $B=\bigcup\limits_{m<\omega}
\lim\big((\name{T}_m)^G\big)$ is  an $\bar{\cO}^6$--{\bf pots} set. Then, by
Proposition \ref{eqnd},  $\NDRK_{\bar{\cO}^6}(\bar{T})=\infty$. Using Lemma 
\ref{lemonrk}(5), by induction on $j<\omega$  we choose $\bm_j, \bm_j^*\in
\bM_{\bar{T},\bar{\cO}^6}$ and $p_j\in G$ such that   
\begin{enumerate}
\item[(i)] $\ndrk_{\bar{\cO}^6}(\bm_j)\geq\omega_1$, $|u^{\bm_j}|>5$ and
  $\bm_j\sqsubset \bm^*_j \sqsubset \bm_{j+1}$,  
\item[(ii)] for each $\nu\in u^{\bm^*_j}$ the set $\{\eta\in u^{\bm_{j+1}}: 
  \nu\vtl \eta\}$ has at least two elements, and
\item[(iii)] $p_j\leq p_{j+1}$, $\ell^{\bm_j}<\ell^{\bm^*_j}= n^{p_j}<
  \ell^{\bm_{j+1}}$ and  $\rng(h_i^{\bm_j})\subseteq M^{p_j}$ for all
  $i<6$,  and   
\item[(iv)] $|\{\eta\rest n^{p_j}:\eta\in u^{\bm_{j+1}}\}|=|u^{\bm_j}|=|
  u^{\bm^*_j}|$.  
\end{enumerate}
To carry out the construction we proceed as follows. Suppose we have
determined $\bm_j$ so that $\ndrk_{\bar{\cO}^6}(\bm_j)\geq \omega_1$. Using
densities given in Claim \ref{cl8}, we find $p_j\in G$ with
$n^{p_j}>\ell^{\bm_j}$ and $\rng(h_i^{\bm_j})\subseteq M^{p_j}$ (for
$i<6$). Next we choose $\bn$ such that $\bm_j\sqsubset \bn$,
$\ndrk_{\bar{\cO}^6} (\bn)\geq \omega_1$, and $\ell^{\bn}>n^{p_j}$. Using
Lemma  \ref{lemonrk}(8) (for a $u'\subseteq u^{\bn}$ such that $\{\eta\rest
\ell^{\bm_j}:\eta\in u'\} =u^{\bm_j}$, $|u'|=|u^{\bm_j}|$) we may
also demand that $|u^{\bn}|=|u^{\bm_j}|$. Now we let
\begin{itemize}
\item $\ell=n^{p_j}$, $u=\{\eta\rest \ell: \eta\in u^{\bn}\}$, 
\item $\bar{h}=\langle h_i:i<6\rangle$, where for $i<6$ and  $(\eta,\nu)\in
  \big(u^{\bn}\big)^{\langle 2\rangle}$\\ 
$h_i(\eta\rest\ell,\nu\rest\ell)= h^{\bn}_i(\eta,\nu)=h^{\bm_j}_i(
\eta\rest\ell^{\bm_j},\nu\rest \ell^{\bm_j})$,
\item $\bar{g}=\langle g_i:i<6\rangle$, where for $i<6$ and 
$(\eta,\nu)\in \big(u^{\bn}\big)^{\langle 2\rangle}$\\
$g_i(\eta\rest\ell,\nu\rest\ell)= \big\{\rho\rest\ell:\rho\in
g^{\bn}_i(\eta,\nu)\big\}$.  
\end{itemize}
Clearly, $\bm^*_j=(\ell,6, u,\bar{h},\bar{g})\in \bM_{\bar{T},\bar{\cO}^6}$
and $\bm_j\sqsubset \bm^*_j$. Finally use Lemma \ref{lemonrk}(5) to pick
$\bm_{j+1}\sqsupset \bn$ such that $\ndrk(\bm_{j+1})\geq \omega_1$ and
condition (ii) is satisfied. Note that $\bm^*_j\sqsubset \bm_{j+1}$. 
\medskip

Then, by (iii)+(iv), $\bm_j,\bm^*_j\in {\mathbf
  M}^{n^{p_j}}_{\bar{t}^{p_j},\bar{\cO}^6}$. It follows from Claim \ref{cl7}
that for some $w_j\subseteq w^{p_j}$ and $\rho_j\in {}^{n^{p_j}}2$ we have 
$(n^{p_j},w_j,\bm_j^*+\rho_j) \in \cM^{p_j}$.  
\medskip

Fix $j$ for a moment and consider $(n^{p_j}, w_j,\bm^*_j+\rho_j ) \in
\cM^{p_j}\subseteq \cM^{p_{j+1}}$ and $(n^{p_{j+1}}, w_{j+1},\bm^*_{j+1}
+\rho_{j+1} )\in \cM^{p_{j+1}}$. (Note that since $(n^{p_j},w_j,\bm^*_j+
\rho_j)\in \cM^{p_j}$, we know that $r_{h^{\bm^*_j}_i(\eta,\nu)}\leq n^{p_j}$ for
  all $i<6$, $(\eta,\nu)\in u^{\bm^*_j}$.) Since $(\bm_j^*+(\rho_{j+1}\rest
n^{p_j}))\sqsubset (\bm^*_{j+1}+\rho_{j+1})$, we  may choose $w^*_j\subseteq
w_{j+1}$ such that $(n^{p_j},w^*_j,\bm_j^*+(\rho_{j+1}\rest n^{p_j}))\in
\cM^{p_{j+1}}$.  Since $(\bm^*_j+\rho_j)+ (\rho_j+\rho_{j+1}\rest
n^{p_j})=\bm_j^*+(\rho_{j+1}\rest n^{p_j})$, we may use clause $(*)_{10}$
for $p_{j+1}$ to conclude that  $\rk(w^*_j)=\rk(w_j)$.

Condition (ii) of the choice of $\bm_{j+1}$ implies that 
\[(\forall \gamma\in w^*_j)(\exists\delta\in w_{j+1}\setminus
  w_j^*)( \eta^{p_{j+1}}_\gamma\rest n ^{p_j}= \eta_\delta^{p_{j+1}}
  \rest n^{p_j}).\] 
Let $\delta(\gamma)$ be the smallest $\delta\in w_{j+1}\setminus w_j^*$
with the above property and let $w^*_j(\gamma)=(w^*_j\setminus
\{\gamma\})\cup \{\delta(\gamma)\}$. Then, for $\gamma\in w^*_j$,
$(n^{p_j}, w^*_j(\gamma), \bm_j^*+(\rho_{j+1}\rest n^{p_j}))\in
\cM^{p_{j+1}}$ and therefore, by clause $(*)_{10}$ for $p_{j+1}$, we get
that for each $\gamma\in w_j$: 
\[\rk(w^*_j(\gamma))=\rk(w^*_j),\quad \zeta(w^*_j(\gamma))
  =\zeta(w_j^*), \quad\mbox{ and }\quad k(w^*_j(\gamma))
  =k(w_j^*).\]
Let $n=|w^*_j|$, $\zeta=\zeta(w_j^*)$, $k=k(w_j^*)$, and let
$w_j^*=\{\alpha_0,\ldots,\alpha_k,\ldots, \alpha_{n-1}\}$ be the increasing
enumeration. Let $\alpha^*_k=\delta(\alpha_k)$. Then clause $(*)_{10}$ also 
gives that $w_j^*(\alpha_k)=\{\alpha_0,\ldots,\alpha_{k-1},\alpha_k^*,
\alpha_{k+1}, \ldots, \alpha_{n-1}\}$ is the increasing enumeration.  Now, 
\[\begin{array}{l}
\bbM\models R_{n,\zeta}[\alpha_0,\ldots,\alpha_{k-1},\alpha_k,
\alpha_{k+1}, \ldots, \alpha_{n-1}]\qquad \mbox{ and}\\
\bbM\models R_{n,\zeta}[\alpha_0,\ldots,\alpha_{k-1},\alpha_k^*,
\alpha_{k+1}, \ldots, \alpha_{n-1}],
\end{array}\]
and consequently if $\rk(w^*_j)\geq 0$, then 
\[\rk(w_{j+1})\leq \rk(w^*_j\cup\{\alpha^*_k\})<\rk(w^*_j)=\rk(w_j)\] 
(remember $(\circledast)_{\rm d}$ from the very beginning of the proof of
the Theorem). 
\medskip

Now, unfixing $j$, it follows from the above considerations that for some 
$j_0<\omega$ we must have:
\begin{enumerate}
\item[(a)] $\rk(w^*_{j_0})=-1$, and  
\item[(b)] $(n^{p_{j_0}}, w^*_{j_0}, \bm_{j_0}^*+(\rho_{j_0+1}\rest
  n^{p_{j_0}})), (n^{p_{j_0+1}}, w_{j_0+1},\bm^*_{j_0+1}+\rho_{j_0+1}) \in
  \cM^{p_{j_0+1}}$,  
\item[(c)] for each $\nu\in u^{\bm^*_{j_0}}$ the set $\{\eta\in u^{\bm^*_{j_0+1}}: 
  \nu\vtl \eta\}$ has at least two elements.
\end{enumerate}
However, this contradicts clause $(*)_{11}$ (for $p_{j_0+1}$).
\end{proof}
\end{proof}

\section{Conclusions and Questions}

\begin{corollary}
\label{corforotps}
Assume ${\rm NPr}_{\omega_1}(\lambda)$ and $\lambda=\lambda^{\aleph_0}
<\mu= \mu^{\aleph_0}$.
\begin{enumerate}
\item Let $\bar{\cO}$ be a nice indexed base. Then there 
is a ccc forcing notion $\bbQ$ of size $\mu$ forcing that:
\begin{itemize}
\item $2^{\aleph_0}=\mu$ and there is a $\Sigma^0_2$ set $B\subseteq\can$
  which has $\lambda$ many pairwise $\bar{\cO}$--nondisjoint translates but
  does not have  $\lambda^+$ many pairwise $\bar{\cO}^6$--nondisjoint
  translates.   
\end{itemize}
\item In particular, there is a ccc forcing notion $\bbQ'$ of size $\mu$
  forcing that:   
\begin{itemize}
\item $2^{\aleph_0}=\mu$ and for some $\Sigma^0_2$ set
  $B\subseteq\can$ there are pairwise distinct $\langle
  \eta_\xi:\xi<\lambda\rangle$ such that $(B+\eta_\xi)\cap (B+\eta_\zeta)$
  is uncountable for each $\xi,\zeta<\lambda$, but
\item for any set $A\subseteq \can$ of size $\lambda^+$ there are $x,y\in A$
  such that $|(B+x)\cap (B+y)|<6$. 
\end{itemize}
\end{enumerate}
\end{corollary}

\begin{proof}
(1) Let $\bbP$ be the forcing notion given by Theorem \ref{522fortranslate}
and let $\bbQ=\bbP*\bbC_\mu$. The set $B$ added by $\bbP$ is a
$\bar{\cO}^6$--{\bf npots}--set in $\bV^\bbP$, so by Proposition \ref{eqnd}
we got $\NDRK_{\bar{\cO}^6}(\bar{T})=\infty$. The rank
$\ndrk_{\bar{\cO}^6}^{\bar{T}}$ is absolute, so in $\bV^\bbQ$ we still have 
$\NDRK_{\bar{\cO}^6}(\bar{T})=\infty$ and thus $B$ is a $\bar{\cO}^6$--{\bf
  npots}--set in $\bV^\bbQ$. By  \ref{proptostart}(3) this set  cannot have
$\lambda^+$ pairwise $\bar{\cO}^6$--nondisjoint translates, but it does have
$\lambda$ many pairwise $\bar{\cO}$--nondisjoint translates (by
absoluteness).  
\end{proof}

\begin{corollary}
\label{MAgives}
Assume ${\bf MA}$ and $\aleph_\alpha<\con$, $\alpha<\omega_1$. 
\begin{enumerate}
\item Let $\bar{\cO}$ be a nice indexed base. Then there exists a
  $\Sigma^0_2$ $\bar{\cO}^6$--{\bf npots}--set $B\subseteq\can$ which has
  $\aleph_\alpha$  many pairwise  $\bar{\cO}$--nondisjoint translations.
\item In particular, there exists a $\Sigma^0_2$ set $B\subseteq \can$ such
  that
  \begin{itemize}
  \item for some pairwise distinct $\langle \eta_\xi:\xi<\aleph_\alpha
    \rangle \subseteq \can$ the intersections $(B+\eta_\xi)\cap
    (B+\eta_\zeta)$ are uncountable for each $\xi,\zeta<\aleph_\alpha$, but
\item for every perfect set $P\subseteq \can$ there are $x,y\in P$ such that
  $|(B+x)\cap (B+y)|<6$.    
  \end{itemize}
\end{enumerate}
\end{corollary}

\begin{proof}
Standard consequence of the proof of Theorem \ref{522fortranslate}, using
the fact that ``$B$ is a $\bar{\cO}^6$--{\bf npots}--set'' is sufficiently
absolute by Proposition \ref{eqnd}.
\end{proof}

\begin{problem}
  \begin{enumerate}
\item Can one differentiate between various nice $\bar{\cO}$ in the context
  of our results? In particular:
  \item Is it consistent that for some nice $\bar{\cO}$  there is an $\Sigma^0_2$
  $\bar{\cO}$--{\bf npots}--set which has   $\aleph_\alpha$  many pairwise
  $\bar{\cO}$--nondisjoint translations, but for some other nice
  $\bar{\cO}^*$ every $\Sigma^0_2$ set with $\aleph_\alpha$  many pairwise
  $\bar{\cO}^*$--nondisjoint translations is automatically $\bar{\cO}^*$--{\bf pots} ?
\item Is it consistent that there is an $\Sigma^0_2$ set  $B\subseteq\can$
  which is has  $\aleph_\alpha$  many pairwise  $\bar{\cO}^\per$--nondisjoint 
  translations, is $\bar{\cO}^\per$--{\bf npots}, but is also
  $\bar{\cO}^6$--{\bf pots}?
  \end{enumerate}
\end{problem}

\begin{problem}
  \begin{enumerate}
  \item Consider the forcing notion $\bbP$ given by Theorem
    \ref{522fortranslate} for $\bar{\cO}^\per$. In the forcing extension by
    $\bbP$, the ranks $\NDRK_{\bar{\cO}^6}(\bar{T})$ and
    $\NDRK_{\bar{\cO}^\per}(\bar{T})$ are both countable. Are they equal? What
    are their values?   
  \item Does there exist a sequence of trees $\bar{T}^*$ (as in Assumptions
    \ref{hyp1}) for which  the ranks $\NDRK_{\bar{\cO}^\per}(\bar{T})$ and
    $\NDRK_{\bar{\cO}^\iota}(\bar{T})$ are different (for some/all $\iota$)?
\item Generalize the construction of \cite{RoSh:1170} to arbitrary nice
  $\bar{\cO}$.
\item Generalize the result of the present paper to the context of arbitrary
  perfect Abelian Polish groups.
  \end{enumerate}
\end{problem}


\end{document}